\documentclass[twocolumn]{autart}    

\usepackage{graphicx}          
\usepackage[utf8]{inputenc}
\usepackage{textcomp}

\usepackage{graphicx}
\usepackage{amsmath}

\usepackage{siunitx}
\usepackage{longtable,tabularx}
\usepackage{subfigure}
\usepackage{algorithm}
\usepackage{algpseudocode}
\usepackage{xcolor}
\usepackage{comment}
\usepackage{upgreek}
\usepackage{kotex}
\usepackage{graphicx}
\usepackage{booktabs}
\usepackage{hyperref}
\usepackage{setspace}
\usepackage{mathtools}
\usepackage{amsfonts}
\usepackage{natbib}
\usepackage{enumitem}

\graphicspath{{images/}}

\begin{document}

\begin{frontmatter}

\title{On the Reachability of 3-Dimensional Paths with a Prescribed Curvature Bound\thanksref{footnoteinfo}} 

\thanks[footnoteinfo]{This paper was not presented at any conference.}

\author[KAIST]{Juho Bae}\ead{johnbae1901@kaist.ac.kr}, 
\author[Nearthlab]{Ji Hoon Bai}\ead{jihoon.bai@nearthlab.com}, 
\author[Nearthlab]{Byung-Yoon Lee}\ead{byungyoon.lee@nearthlab.com}, 
\author[Nearthlab]{Jun-Yong Lee}\ead{junyong.lee@nearthlab.com}, 
\author[KAIST]{Chang-Hun Lee\thanksref{corrautinfo}}\ead{lckdgns@kaist.ac.kr}

\address[KAIST]{Korea Advanced Institute of Science and Technology, Daejeon, 34141, Republic of Korea}
\address[Nearthlab]{Nearthlab, Seoul, 05836, Republic of Korea}

\thanks[corrautinfo]{Corresponding author.}
          

\begin{abstract}                          
This paper presents the reachability analysis of curves in $\mathbb{R}^3$ with a prescribed curvature bound. Based on Pontryagin Maximum Principle, we leverage the existing knowledge on the structure of solutions to minimum-time problems, or Markov-Dubins problem, to reachability considerations. Based on this development, two types of reachability are discussed. First, we prove that any boundary point of the reachability set, with the directional component taken into account as well as geometric coordinates, can be reached via curves of H, CSC, CCC, or their respective subsegments, where H denotes a helicoidal arc, C a circular arc with maximum curvature, and S a straight segment. Second, we show that the reachability set when directional component is not considered\textemdash{}the position reachability set\textemdash{}is simply a solid of revolution of its two-dimensional counterpart, the Dubins car. These findings extend the developments presented in literature on Dubins car into spatial curves in $\mathbb{R}^3$.
\end{abstract}
\end{frontmatter}
%
\section{Introduction}
Analysis of reachability of dynamical systems is crucial in a wide range of applications~(\cite{matthias2010}). Accordingly, various studies have explored methods for constructing the reachability set of curves in two-dimensional (2D) plane with prescribed curvature bound~(\cite{Melzak1961},~\cite{wong1974},~\cite{cockayne1975},~\cite{patsko2003},~\cite{patsko2022}). The emphasis on the curvature-bounded paths is because the dynamics of planar curves with a prescribed curvature bound, or \textit{Dubins car}~(\cite{dubins1957}), represents numerous practical platforms.

In~\cite{dubins1957}, the author addressed the 2D case of the problem first introduced by A. A. Markov, also known as \textit{Markov-Dubins problem}, which questions the structure of the shortest paths in $\mathbb{R}^n$ with prescribed curvature bounds along with specified initial and terminal locations and directions. Dubins proved that the minimizers are curves of CSC, CCC, or their subsegments in 2D, where C denotes a circular arc with maximum curvature and S a straight segment. (Throughout this paper, we will refer to the collection of all CSC curves as \textit{class CSC}. This notation will be applied strictly in a sense that, for instance, any CSC curve with the S component having zero length will be referred to as CC, not CSC. Subsegments of, for instance, class CSC, refer to classes C, S, CS, SC, and CC.) 
While the central arguments in \cite{dubins1957} were based on geometric approaches, subsequent works have reformulated the problem as an OCP and applied the Pontryagin Maximum Principle (PMP) to reproduce the same results~(\cite{Pecsvaradi1972},~\cite{johnson1974},~\cite{Boissonnat1992}). 

These works facilitated subsequent studies on reachability of Dubins car. The question was as follows~(\cite{cockayne1975}): \textit{``A particle moves in $\mathbb{R}^2$ with constant speed and subject to an upper bound on the curvature of its path. At each time $t$, what is the set of all possible positions for the particle?''} The small-time case was addressed in~\cite{Melzak1961}, followed by~\cite{cockayne1975} for general $t$. 
By applying geometric arguments,~\cite{cockayne1975} first introduced an implicit description of the reachability that every boundary point can be reached by curves in classes CS, CC, and their subsegments. Hence, collecting the endpoints of all such curves provides a small set containing the boundary. Then, an additional condition that the reachability set varies continuously in time facilitated complete and analytical description of the boundary of the reachability set. In this problem, the target of reachability analysis was only the geometric coordinates, which we refer to as \textit{position reachability set}.

Reachability analysis of Dubins car with directional components considered was first addressed in~\cite{patsko2003}. PMP was used to conclude that every boundary point of the reachability set can be reached by curves in classes CSC, CCC, or their subsegments. This implicit description does not solely provide complete description of the boundary of the reachability set for general $t$, but was sufficient for construction for small-time. Subsequently, a follow-up study in~\cite{patsko2022} provided additional conditions to facilitate complete and analytical description of the reachability set for all $t$.

Prior developments on the reachability of the Dubins car~(\cite{cockayne1975},~\cite{patsko2003}) enabled numerous subsequent advancements~(\cite{buzikov2021},~\cite{Buzikov2024},~\cite{buzikov2022},~\cite{chen2023},~\cite{patsko2022}). Motivated by these works, we extend the existing results on Dubins car into 3D curves with prescribed curvature bound. Since names such as \textit{Dubins airplane} or \textit{3D Dubins car} have already been taken in the literature to represent different dynamics\textemdash{}extension of Dubins car into 3D by including an independent altitude component~(\cite{chitsaz2007})\textemdash{}we will refer to our dynamics as \textit{3D curves with prescribed curvature bound}. Unlike `Dubins airplane', our dynamics maintains the curvature bound property. As reachability analysis of Dubins car in literature utilized the developments on its minimum time problem, we similarly leverage the results of minimum time problem for the 3D case studied in~\cite{sussmann1995}.

Contributions of this paper are as follows: First, we reveal the internal connection between reachability and time optimality which allows to apply the conclusions in~\cite{sussmann1995} on the minimum time problem to reachability considerations. Second, we prove that the boundary points of the reachability set (i.e., when directional components are also considered.) can be reached via CCC, CSC, their subsegments, or H, where H denotes certain class of helicoidal arcs. This extends the results in~\cite{patsko2003} into dimension 3. Last, we present complete and analytical description of the reachability set of 3D curves with prescribed curvature bound. We prove that the position reachability set in this case can be simply generated by rotating the position reachability set of Dubins car in~\cite{cockayne1975}. This extends the results made in~\cite{cockayne1975} into dimension 3.
%
\section{Preliminary}\label{sec:02}
\subsection{Definitions and preliminary theorems}
Consider the process in $\mathbb{R}^{n_{\boldsymbol{x}}}$
\begin{equation} \label{eq:dynamics}
	\dot{\boldsymbol{x}} = f(\boldsymbol{x}, \boldsymbol{u}), \quad \boldsymbol{x}(0) = \boldsymbol{x}_0
\end{equation}
with $\boldsymbol{x} \in \mathbb{R}^{n_{\boldsymbol{x}}}$ and $\boldsymbol{u} \in \mathbb{R}^{n_{\boldsymbol{u}}}$, where $f$ is $C^1$ in $\mathbb{R}^{n_{\boldsymbol{x}} + n_{\boldsymbol{u}}}$. $\boldsymbol{x}_0 \in \mathbb{R}^{n_{\boldsymbol{x}}}$ represents the initial state. Let $\mathcal{F}$ be the family of all measurable functions $\boldsymbol{u}$ on a specified interval $[0, T]$ with values $\boldsymbol{u}(t) \in \Omega$ for some compact restraint set in $\mathbb{R}^{n_{\boldsymbol{u}}}$. We introduce the following assumptions on $f$.
\begin{assum} \label{assum:01}
	\hfill
	\begin{itemize}
		\item For each $\boldsymbol{u} \in \mathcal{F}$, there exists a response on $[0, T]$, 
		\item there exists a uniform bound $|\boldsymbol{x}(t)| < m$ for all responses on $[0, T]$, 
		\item $V(\boldsymbol{x}) = \{ f(\boldsymbol{x}, \boldsymbol{u}) : \boldsymbol{u} \in \Omega \}$ is convex for each $\boldsymbol{x}$.
	\end{itemize}
\end{assum}
In this context, \textit{uniform bound} means that the choice of the upper bound $m$ is independent of the choice of controller $\boldsymbol{u} \in \mathcal{F}$ and time $t$. Compactness of $V(\boldsymbol{x})$ follows trivially from compactness of $\Omega$ and continuity of $f$. Assumption~\ref{assum:01} is sufficient for our discussions though some can be relaxed beyond the purpose of applications to 3D curves with prescribed curvature bounds. For $t \in [0, T]$, the \textit{Reachability Set} $\mathcal{G}(t)$ is defined as follows.
\begin{defn}
\begin{equation}
	\mathcal{G}(t) \equiv \{ \boldsymbol{x}(t) : \boldsymbol{u} \in \mathcal{F} \}
\end{equation}
\end{defn}
Assumption~\ref{assum:01} and the well-known Filippov's Theorem assert that each $\mathcal{G}(t)$ is compact. Moreover, under Assumption~\ref{assum:01}, it was proved in~\cite{Lee1967} that $\mathcal{G}(t)$ varies continuously in time on $[0, T]$, where continuity is with respect to Hausdorff metric.
It was further proved in~\cite{Lee1967} the following theorem, where boundary of a set $S$ is denoted as $bd(S)$.
\begin{thm}{(\cite{Lee1967})} \label{thm:1}
	Let $\boldsymbol{u}^*(t) \in \mathcal{F}$ have a response $\boldsymbol{x}^*(t)$ such that $\boldsymbol{x}^*(t_f) \in bd(\mathcal{G}(t_f))$. Then there exists a nontrivial adjoint system $\boldsymbol{p}^*(t)$ such that 
	\begin{equation} \label{eq:thm1-1}
		\dot{\boldsymbol{p}}^*(t) = -\boldsymbol{p}^*(t)^T \frac{\partial f}{\partial \boldsymbol{x}}( \boldsymbol{x}^*(t), \boldsymbol{u}^*(t) ),
	\end{equation}
	\begin{equation} \label{eq:thm1-2}
	\begin{split}
		& \mathcal{H}( \boldsymbol{p}^*(t), \boldsymbol{x}^*(t), \boldsymbol{u}^*(t) ) = \max\limits_{\boldsymbol{u} \in \Omega} \mathcal{H}( \boldsymbol{p}^*(t), \boldsymbol{x}^*(t), \boldsymbol{u} ) \\
		& \text{for a.e. } t 
	\end{split}
	\end{equation}
	where $\mathcal{H}( \boldsymbol{p}, \boldsymbol{x}, \boldsymbol{u} ) \equiv \boldsymbol{p}^T f( \boldsymbol{x}, \boldsymbol{u} )$. 
\end{thm}
The `almost everywhere' condition is abbreviated as \textit{a.e.} $t$. A slightly stronger statement is available if all responses lie on a smooth manifold $M$. In this case, since each $f(\boldsymbol{x}^*, \boldsymbol{u}^*)$ is a tangent vector in $T_{\boldsymbol{x}^*}M$, the adjoint response $\boldsymbol{p}^*$ is viewed as a covector on $M$ rather than a vector in the ambient space $\mathbb{R}^{n_{\boldsymbol{x}}}$. Consequently, the nontriviality condition of $\boldsymbol{p}^*$ is applied by interpreting $\boldsymbol{p}^*$ as a covector on $M$, instead of merely being a nonzero vector in $\mathbb{R}^{n_{\boldsymbol{x}}}$. In other words, $\boldsymbol{p}^*$ has to be considered as a nonvanishing linear functional $\boldsymbol{z} \mapsto \langle \boldsymbol{p}^*, \boldsymbol{z} \rangle$ from $T_{\boldsymbol{x}^*} M$ to $\mathbb{R}$. Benefit of such distinction between the nontriviality conditions will be apparent in the subsequent sections. A more detailed description regarding the formulation of maximum principle under manifold setting can be found in~\cite{sussmann1998}. Under this consideration, the following theorem directly follows from Theorem~8.2 in~\cite{sussmann1998}.
\begin{thm}{(\cite{sussmann1998})} \label{thm:2}
	Suppose the state space is $M = M_1 \times M_2$ so that each response $\boldsymbol{x}(t) = (\boldsymbol{x}_1(t), \boldsymbol{x}_2(t))$ lies on $M_1 \times M_2$, where $M_1$ and $M_2$ are smooth manifolds that $\boldsymbol{x}_1$ and $\boldsymbol{x}_2$ lie respectively. Define a projection map $\pi: \boldsymbol{x} \mapsto \boldsymbol{x}_1$. 
	For a final time $t_f \in [0, T]$, suppose $\boldsymbol{u}^*(t) \in \mathcal{F}$ has a response $\boldsymbol{x}^*(t)$ such that $\pi\left(\boldsymbol{x}^*(t_f)\right) \in bd(\pi\left(\mathcal{G}(t_f)\right))$. Then there exists a nontrivial adjoint system $\boldsymbol{p}^*(t) = (\boldsymbol{p}^*_1(t), \boldsymbol{p}^*_2(t))$ satisfying Eqs.~\eqref{eq:thm1-1} and~\eqref{eq:thm1-2} on $[0, t_f]$ such that $\boldsymbol{p}^*_2(t_f) \in T^*_{\boldsymbol{x}_2^*(t_f)} M_2$ is a zero covector. 
\end{thm}
It is noteworthy that under projection setting, the only additional condition is $\boldsymbol{p}^*_2(t_f)$ being a zero covector. To facilitate further discussions, let us define the set of all solutions $\{ \boldsymbol{x}^*, \boldsymbol{u}^*, \boldsymbol{p}^* \}$ of Eqs.~\eqref{eq:dynamics},~\eqref{eq:thm1-1}, and \eqref{eq:thm1-2} as follows. Alongside, we define the set of endpoints for such solutions. 
\begin{defn}
	Given final time $t_f \in [0, T]$, 
	\begin{enumerate}
		\item $S^*(t_f) \equiv \{ \{ \boldsymbol{x}^*, \boldsymbol{u}^*, \boldsymbol{p}^* \} : \text{nontrivial solutions of } \\ \text{Eqs.~\eqref{eq:dynamics},~\eqref{eq:thm1-1}, and \eqref{eq:thm1-2} on } [0, t_f] \}$
		\item $E^*(t_f) \equiv \{ \boldsymbol{x}^*(t_f) : \{ \boldsymbol{x}^*, \boldsymbol{u}^*, \boldsymbol{p}^* \} \in S^*(t_f) \}$
	\end{enumerate}
\end{defn}
Since $\mathcal{G}(t_f)$ is compact under Assumption~\ref{assum:01}, all of its boundary points are contained in itself and therefore reachable at time $t_f$. Consequently, the following corollary holds for each $t_f$.
\begin{cor} \label{cor:01}
	$bd(\mathcal{G}(t_f)) \subseteq E^*(t_f)$
\end{cor}
Corollary~\ref{cor:01} facilitates the identification of the boundary of $\mathcal{G}(t_f)$ up to certain extent by aggregating all the endpoints of the trajectories in $S^*(t_f)$. However, it should be noted that identifying $E^*(t_f)$ does not completely determine the boundary of $\mathcal{G}(t_f)$ in general, because Theorem~\ref{thm:1} and its byproduct, Corollary~\ref{cor:01}, only serve as necessary conditions but not sufficient. Consequently, some points in $E^*(t_f)$ may lie interior to $\mathcal{G}(t_f)$ and additional information is required to exclude such points to completely determine the boundary. 
 
Additional information contained in Theorem~\ref{thm:1} but not in Corollary~\ref{cor:01} is Theorem~\ref{thm:1} states that for `every' trajectory reaching the boundary, there exists a nontrivial solution of Eqs.~\eqref{eq:thm1-1} and~\eqref{eq:thm1-2}. Hence, the points in $E^*(t_f)$ that can be reached by an alternative controller that is not maximal (i.e., violates the condition in Eq.~\eqref{eq:thm1-2}) lie interior to $\mathcal{G}(t_f)$. This was leveraged in the reachability analysis of Dubins car in~\cite{patsko2003} (e.g., Lemma~2).

Another information, which in this case not contained in Theorem~\ref{thm:1}, is that $\mathcal{G}(t)$ varies continuously in time. This facilitated complete description of the position reachability set of Dubins car for all $t$~(\cite{cockayne1975}). Throughout our reachability analysis in Section~\ref{sec:03}, we will utilize these two additional information extensively.
\subsection{Markov-Dubins problem in dimensions 2 and 3}
Formulation of Markov-Dubins problem is as follows: given a constant $\kappa_M > 0$, two points $\boldsymbol{x}_0$, $\boldsymbol{x}_1 \in \mathbb{R}^n$, and unit vectors $\boldsymbol{y}_0$, $\boldsymbol{y}_1 \in \mathbb{R}^n$, find a curve $\gamma: [0, L] \mapsto \mathbb{R}^n$ of shortest length (i.e., minimum $L$) such that (1) $\gamma(s)$ is parametrized by arc length, (2) $\gamma(0) = \boldsymbol{x}_0$, $\gamma'(0) = \boldsymbol{y}_0$, $\gamma(L) = \boldsymbol{x}_1$, and $\gamma'(L) = \boldsymbol{y}_1$, and (3) $\gamma'$ is absolutely continuous and satisfies the curvature bound: $|\gamma''(s)| \leq \kappa_M$ for almost all $s \in [0, L]$.

The case of $n = 2$, also known as \textit{Dubins car}, was first studied in~\cite{dubins1957} by geometric means proving that the minimizers are of CSC, CCC, and their subsegments. After~\cite{dubins1957}, the problem was formulated as a minimum time optimal control problem and treated by means of PMP in~\cite{Pecsvaradi1972},~\cite{johnson1974}, and~\cite{Boissonnat1992} to reproduce the same results. The problem setup was as follows, taking $\kappa_M = 1$ for simplicity.

Minimize
\begin{equation}
	\enspace J  =  \int_{0}^{t_f} 1 dt 
\end{equation}
subject to
\begin{equation} \label{eq:2d_dynamics}
	\dot{\boldsymbol{\chi}}\left( t \right) = \begin{pmatrix} \cos\theta(t) \\ \sin\theta(t) \\ 0 \end{pmatrix} + \begin{pmatrix} 0 \\ 0 \\ 1 \end{pmatrix} u(t), \quad |u| \leq 1
\end{equation}
\begin{equation} 
	\begin{split}
		& \boldsymbol{\chi}( 0 ) = \boldsymbol{\chi}_0, \quad \boldsymbol{\chi}( t_f ) = \boldsymbol{\chi}_1, \\
		& \text{final time } t_f \text{ is free}
	\end{split}
\end{equation}
where the state variable is $\boldsymbol{\chi} = [x^1,\; x^2,\; \theta]^T \in \mathbb{R}^3$ and control variable is $u \in \mathbb{R}$. $x^1$ and $x^2$ represent the geometric coordinates and $\theta$ represents the heading angle. Admissible controls are all measurable functions $t \mapsto u(t) \in [-1, 1]$.

For the case of $n = 3$,~\cite{sussmann1995} provides a detailed description of the structure of the minimum time trajectories through an OCP formulation. Again, by taking $\kappa_M = 1$, the formulation was as follows: 

Minimize
\begin{equation} \label{eq:3d_cost}
	\enspace J  =  \int_{0}^{t_f} 1 dt
\end{equation}
subject to 
\begin{equation} \label{eq:3d_dynamics}
\dot{\boldsymbol{x}} = \boldsymbol{y}, \enspace \dot{\boldsymbol{y}} = \boldsymbol{y} \times \boldsymbol{w}, \enspace \boldsymbol{w} \in \mathbb{B}^3
\end{equation}
\begin{equation} \label{eq:3d_constraints}
	\begin{split}
		& \boldsymbol{x}(0) = \boldsymbol{x}_0, \enspace \boldsymbol{y}(0) = \boldsymbol{y}_0, \\
		& \boldsymbol{x}(t_f) = \boldsymbol{x}_1, \enspace \boldsymbol{y}(t_f) = \boldsymbol{y}_1, \\
		& \text{final time } t_f \text{ is free}
	\end{split}
\end{equation}
where $\boldsymbol{x}$ denotes the geometric coordinates and $\boldsymbol{y}$ the velocity. The control variable, represented by $\boldsymbol{w}$, is restricted to take values in $\mathbb{B}^3$, the closed unit ball in $\mathbb{R}^3$. Admissible controls are all measurable functions $t \mapsto \boldsymbol{w}(t) \in \mathbb{B}^3$. Since we assume parametrization by arc length, it is evident that $\boldsymbol{y}$ must lie on $\mathbb{S}^2$, the unit sphere in $\mathbb{R}^3$. Therefore, the state variable $(\boldsymbol{x}, \boldsymbol{y})$ must lie on the manifold $M = \mathbb{R}^3 \times \mathbb{S}^2$. This formulation on manifold plays a crucial role in the application of nontriviality condition of PMP, but we defer the detailed description until the proof of our main theorem. 

Under this formulation, primarily two arguments were applied to analyze the structure of the minimum time trajectories in~\cite{sussmann1995}: first, solving the PMP conditions of minimum time optimal control problem defined by Eqs.~\eqref{eq:3d_cost},~\eqref{eq:3d_dynamics}, and~\eqref{eq:3d_constraints}, and second, a relatively trivial argument that any trajectory containing a full cycle cannot be a minimum time trajectory. This is because any cycle can be detached from the original trajectory to form a new trajectory with shorter length and unperturbed boundary conditions. We will further refer to this argument as \textit{cycle argument}. Additionally, let us denote by H the collection of all helicoidal arcs with constant curvature 1 and nonvanishing torsion $\tau$ satisfying the ODE in Eq.~\eqref{eq:torsionODE} for some constant $\zeta \in \mathbb{R}$. 
\begin{equation} \label{eq:torsionODE}
	\ddot{\tau} = \frac{3\dot{\tau}^{2}}{2\tau} - 2\tau^3 + 2\tau - \zeta \tau \sqrt{|\tau|}
\end{equation}
The conclusion of~\cite{sussmann1995} was that the minimizers are of classes CSC, CCC, their subsegments, or H with $\zeta \geq 0$.
\subsection{Reachability of Dubins car}
Apart from the Markov-Dubins problem, there had been several attempts to obtain the reachability set of planar curves with prescribed curvature bound. To facilitate subsequent discussions, we introduce the following notations. For each time $t$, we will denote the reachability set of the dynamical system Eq.~\eqref{eq:2d_dynamics} as $\mathcal{G}_2(t)$. By defining a projection map $\pi_2: (x^1, x^2, \theta) \mapsto (x^1, x^2)$, the position reachability set is indicated as $\pi_2\left(\mathcal{G}_2(t)\right)$. Similarly, we denote the reachability set of the dynamical system Eq.~\eqref{eq:3d_dynamics} as $\mathcal{G}_3(t)$. By defining a projection map $\pi_3: (\boldsymbol{x}, \boldsymbol{y}) \mapsto \boldsymbol{x}$, the position reachability set is indicated as $\pi_3\left(\mathcal{G}_3(t)\right)$.

\begin{figure}[htbp]
	\centering
	\subfigure[$t_f = \pi$]{\centering \includegraphics[scale=0.2]{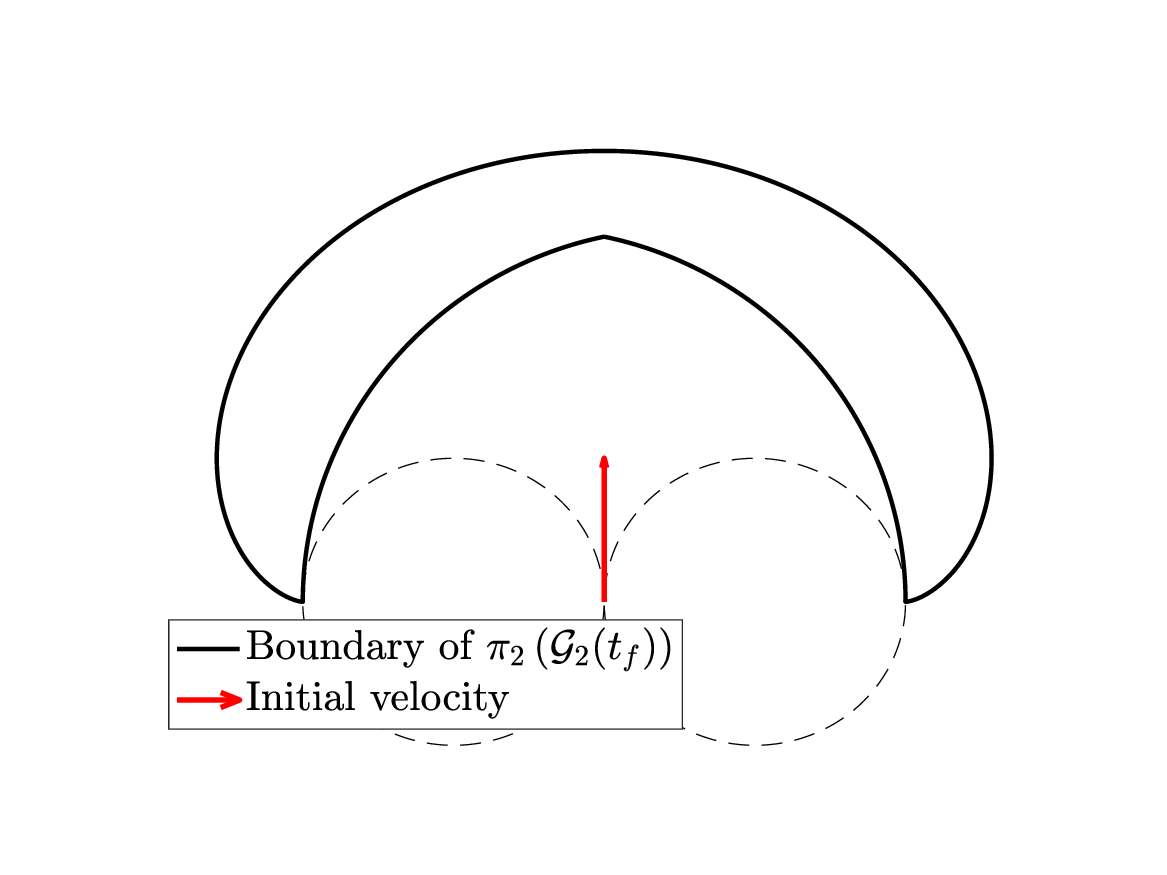}}
	\subfigure[$t_f = 1.5\pi$]{\centering \includegraphics[scale=0.2]{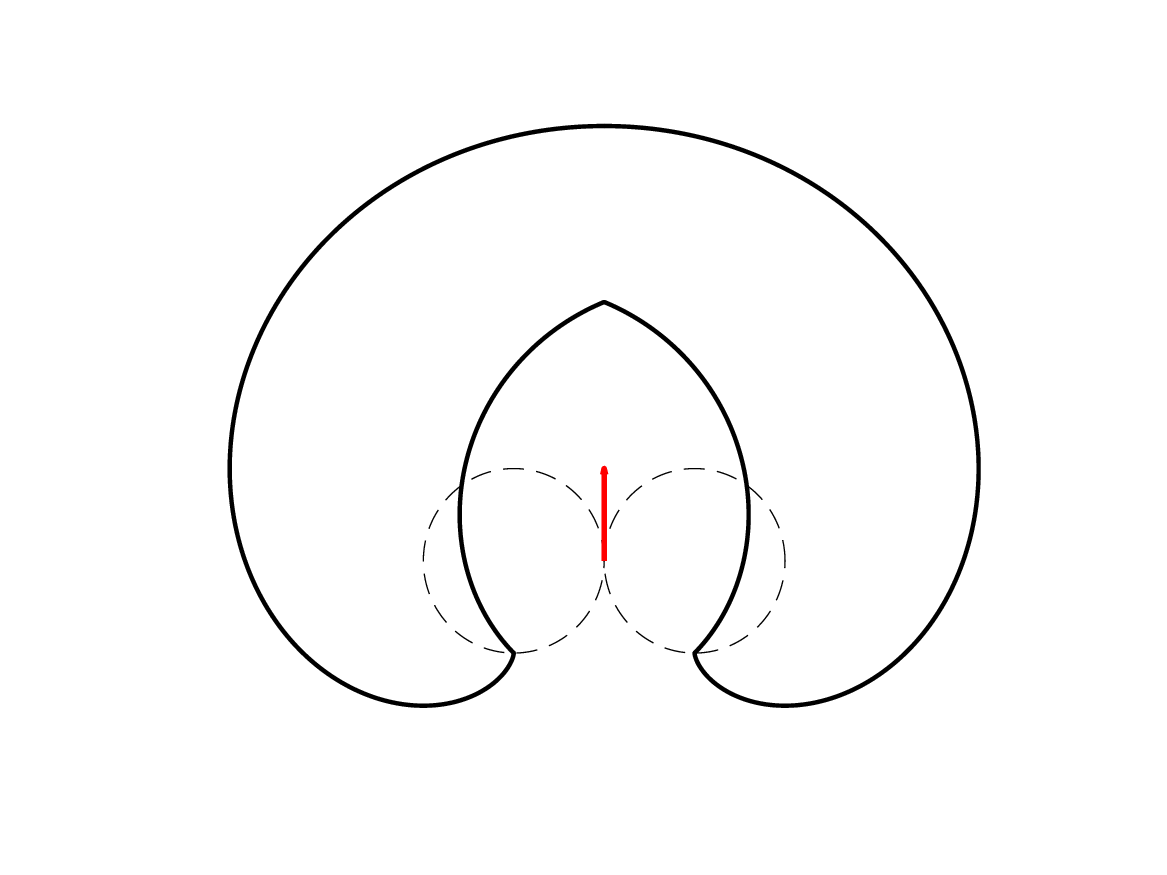}}
	\subfigure[$t_f = 2\pi$]{\centering \includegraphics[scale=0.2]{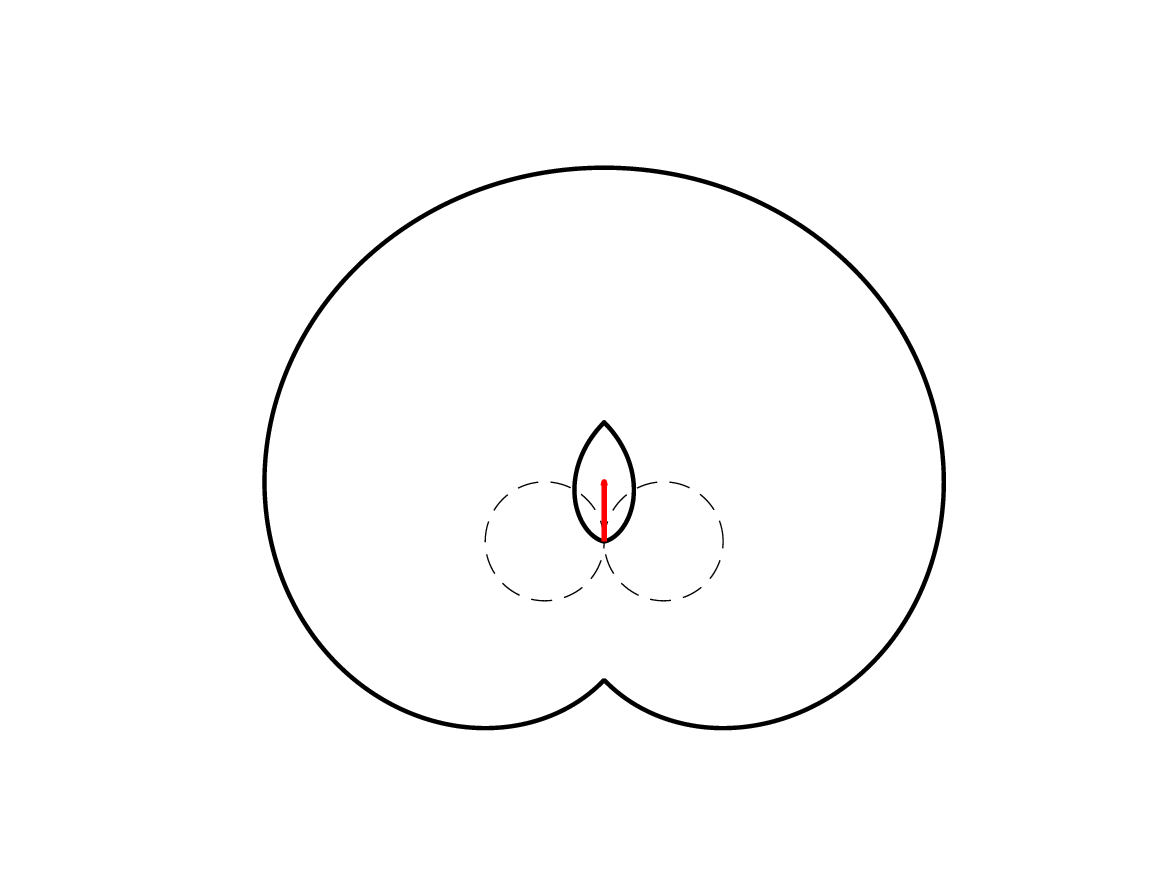}}
	\subfigure[$t_f = 3\pi$]{\centering \includegraphics[scale=0.2]{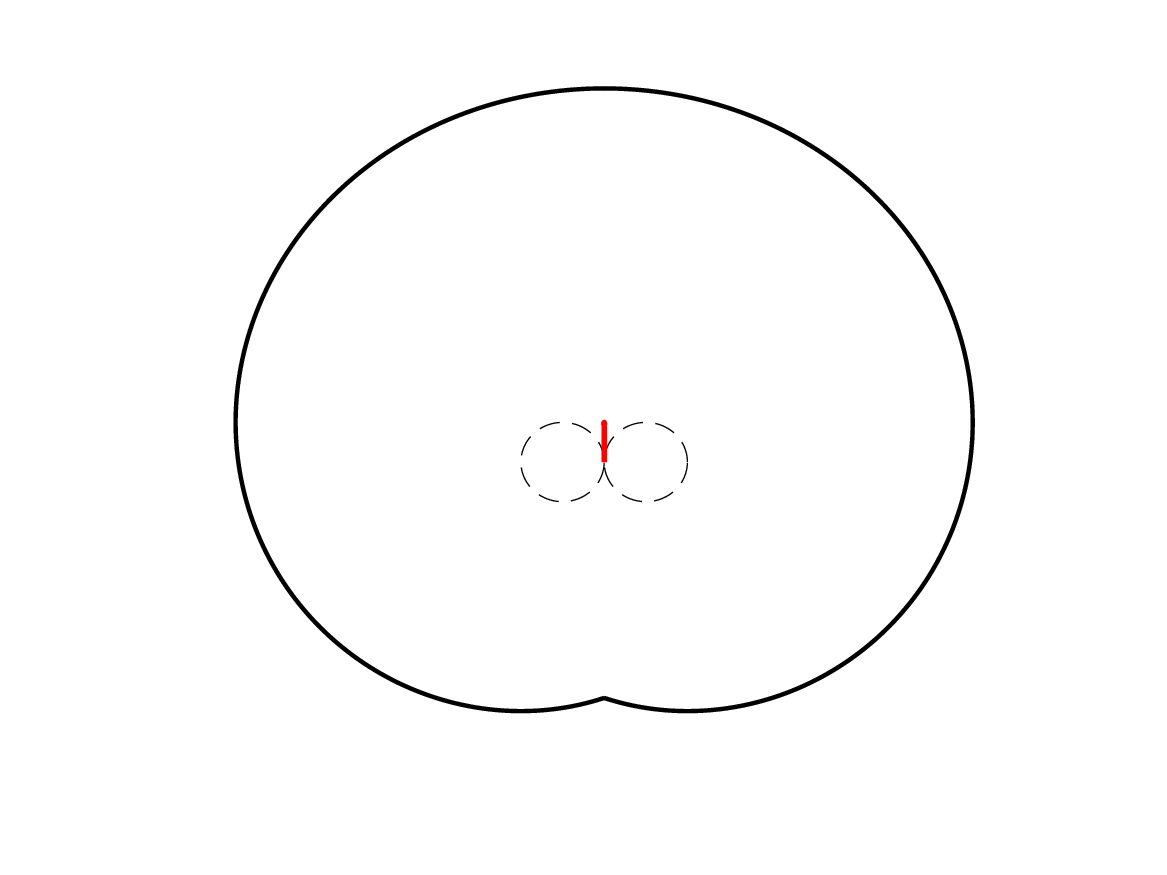}}
	\caption{Evolution of $\pi_2\left(\mathcal{G}_2(t_f)\right)$ over multiple values of $t_f$.}
	\label{fig00}
\end{figure}

The question in~\cite{cockayne1975} can be encoded in terms of our notations as: \textit{``Identify $\pi_2\left(\mathcal{G}_2(t)\right)$.''} It was proved in~\cite{cockayne1975} that every boundary point of $\pi_2\left(\mathcal{G}_2(t)\right)$ can be reached by trajectories of CS, CC, or their subsegments. Identifying a small set containing the boundary points\textemdash{}endpoints of all trajectories in classes CS, CC, and their subsegments\textemdash{}gives a finite number of regions enclosed by arcs generated by such endpoints. However, this in principle does not tell which region belongs to $\pi_2\left(\mathcal{G}_2(t)\right)$ or not because endpoints of some trajectories in such classes may lie interior to $\pi_2\left(\mathcal{G}_2(t)\right)$. For complete determination of the boundary of $\pi_2\left(\mathcal{G}_2(t)\right)$ in~\cite{cockayne1975}, the authors first identified the endpoints of certain trajectories that lie interior to $\pi_2\left(\mathcal{G}_2(t)\right)$. 
By representing the line containing the initial velocity vector by $\ell$, the conclusion can be summarized as follows. 
\begin{prop}{(\cite{cockayne1975})} \label{prop:1}
	Any boundary point of $\pi_2\left(\mathcal{G}_2(t_f)\right)$ can be reached by the curves with length $t_f$ belonging to the following classes: CS, CC, or their subsegments. Moreover, curves of CC required to reach the boundary points are limited to the ones such that 
	\begin{enumerate}
		\item the first C component has length $\leq \frac{\pi}{2}$, 
		\item the second C component has length $\leq 2\pi$, 
		\item the endpoint and the center of the first C component do not lie on the same side of $\ell$. 
	\end{enumerate}
	Similarly, CS curves are limited to the  ones such that 
	\begin{enumerate}
		\item the C component has length $\leq 2\pi$, 
		\item the endpoint and the center of the C component do not lie on the opposite sides of $\ell$. 
	\end{enumerate}
\end{prop}
Let us denote by $E_1(t_f)$ the set of endpoints of curves of length $t_f$ belonging to the classes CS, CC, and their subsegments. The authors first proved that $E_1(t_f)$ contains the boundary of $\pi_2\left(\mathcal{G}_2(t_f)\right)$. Then, it was proved in Theorem 4 in~\cite{cockayne1975} that the points in $E_1(t_f)$ reached by the curves violating the five additional items lie interior to $\pi_2\left(\mathcal{G}_2(t_f)\right)$, thereby proving Proposition~\ref{prop:1}. By denoting these points as $E'(t_f)$, the statement is strengthened as follows: \textit{$E_1(t_f) \setminus E'(t_f)$ contains the boundary of $\pi_2\left(\mathcal{G}_2(t_f)\right)$}. Subsequently, it was confirmed that this leaves finite number of regions enclosed by $E_1(t_f) \setminus E'(t_f)$, which there is a unique way to select the enclosed regions so that their union varies continuously in time. This facilitated the complete determination of $\pi_2\left( \mathcal{G}_2(t_f) \right)$, as depicted in Fig.~\ref{fig00}.
Unlike~\cite{cockayne1975} which relied on geometric arguments, ~\cite{patsko2003} leveraged Theorem~\ref{thm:1} to analyze the reachability set with consideration of directional components, or $\mathcal{G}_2(t)$. It was proved that every boundary point of $\mathcal{G}_2(t)$ can be reached by trajectories of CSC, CCC, or their subsegments. Although Theorem~\ref{thm:1} in principle does not provide complete determination of the boundary, within the specific dynamics of the Dubins car in Eq.~\eqref{eq:2d_dynamics} and small-time $\leq \frac{3.63\pi}{\kappa_M}$, Theorem~\ref{thm:1} was sufficient for construction of the boundary of $\mathcal{G}_2(t)$. This was because for small times, the surfaces generated by collecting the endpoints of the trajectories of CSC, CCC, and their subsegments altogether generated a single surface enclosing a simply-connected set in $\mathbb{R}^3$. (Figs. 8 and 9 in~\cite{patsko2003}) The complete determination of $\mathcal{G}_2(t)$ for general $t$ was finally achieved in the follow-up study of~\cite{patsko2022}.
%
%
\section{Main results}\label{sec:03}
\subsection{Connection with time optimal problems}
In this subsection, we present how the results in~\cite{sussmann1995} can be related with reachability analysis. Let us consider the following OCP which final time is considered as a cost, denoted as $\textbf{P}_{\boldsymbol{x}_f, \phi}$.

[$\textbf{P}_{\boldsymbol{x}_f, \phi}$] : Maximize
\begin{equation} 
	\enspace J  =  \int_{0}^{t_f} \phi dt
\end{equation}
subject to
\begin{equation} 
\begin{split}
	& \dot{\boldsymbol{x}}\left( t \right) = f( \boldsymbol{x}(t), \boldsymbol{u}(t) ), \quad \boldsymbol{u} \in \mathcal{F} \\
	& \boldsymbol{x}( 0 ) = \boldsymbol{x}_0, \quad \boldsymbol{x}( t_f ) = \boldsymbol{x}_f, \\
	& \text{final time } t_f \text{ is free}
\end{split}
\end{equation}
The integrand $\phi \neq 0 \in \mathbb{R}$ is set to be positive for maximum time and negative for minimum time problems. Then, the PMP conditions of $\textbf{P}_{\boldsymbol{x}_f, \phi}$ are as follows: \\ 

A pair $(\boldsymbol{x}^*, \boldsymbol{u}^*)$ defined on $[0, t_f^*] \subseteq [0, T]$ satisfies the PMP if $\exists \{ p^*_{0_{\geq 0}}, \boldsymbol{p}^* \}$, nontrivial, such that 
\begin{enumerate}
	\item Costate differential: $\dot{\boldsymbol{p}}^*(t) = -\boldsymbol{p}^*(t)^T \frac{\partial f}{\partial \boldsymbol{x}}( \boldsymbol{x}^*(t), \boldsymbol{u}^*(t) )$
	\item Pointwise maximum: $\mathcal{H}( p^*_0, \boldsymbol{p}^*(t), \boldsymbol{x}^*(t), \boldsymbol{u}^*(t) ) = \max\limits_{\boldsymbol{u} \in \Omega} \mathcal{H}( p^*_0, \boldsymbol{p}^*(t), \boldsymbol{x}^*(t), \boldsymbol{u} )$ for a.e. $t$
	\item $\mathcal{H}( p^*_0, \boldsymbol{p}^*(t), \boldsymbol{x}^*(t), \boldsymbol{u}^*(t) ) = 0$ for $\forall t \in [0, t^*_f]$
\end{enumerate}
where Hamiltonian is defined as $\mathcal{H}( p_0, \boldsymbol{p}, \boldsymbol{x}, \boldsymbol{u} ) \equiv \boldsymbol{p}^T f( \boldsymbol{x}, \boldsymbol{u} ) + p_0 \phi$. A pair $\{ p^*_{0}, \boldsymbol{p}^* \}$ is said to be nontrivial if $p^*_{0} \neq 0$ or $\boldsymbol{p}^*$ is nonvanishing as a covector on $M$. Any response satisfying the PMP conditions on an interval $[0, t_f^*]$ will be referred to as \textit{maximal response} and the corresponding controller as \textit{maximal controller}. The triplet $\{ \boldsymbol{x}^*, \boldsymbol{u}^*, \boldsymbol{p}^* \}$ is referred to as \textit{maximal}.

It is easy to notice that the costate differential and pointwise maximum conditions of $\textbf{P}_{\boldsymbol{x}_f, \phi}$ problems are equivalent with Eqs.~\eqref{eq:thm1-1} and~\eqref{eq:thm1-2}. Moreover, from everywhere vanishing Hamiltonian condition, vanishing $\boldsymbol{p}^*$ as a covector of $M$ implies $p_0^* = 0$. Thus, $\boldsymbol{p}^*$ has to be nontrivial by itself. This implies that a maximal of a $\textbf{P}_{\boldsymbol{x}_f, \phi}$ problem becomes a nontrivial solution to Eqs.~\eqref{eq:thm1-1} and~\eqref{eq:thm1-2}. Conversely, $\boldsymbol{p}^*(t)^T f( \boldsymbol{x}^*(t), \boldsymbol{u}^*(t) )$ remains constant over time since $f$ is an autonomous system, and existence of $p_0^*$ satisfying the everywhere vanishing Hamiltonian condition follows thereby. Consequently, collecting all nontrivial solutions of Eqs.~\eqref{eq:thm1-1} and~\eqref{eq:thm1-2} is equivalent with finding all maximals of $\textbf{P}_{\boldsymbol{x}_f, \phi}$ problems. In other words, identifying the maximals of $\textbf{P}_{\boldsymbol{x}_f, \phi}$ problems provides every information to construct the set $E^*(t_f)$. This enables to utilize the existing results on time optimal problems to obtain useful conclusions on reachability sets. \cite{sussmann1995} uses PMP to find the maximals of minimum time problem (i.e., $\textbf{P}_{\boldsymbol{x}_f, \phi}$ with $\phi < 0$ described by Eqs.~\eqref{eq:3d_cost},~\eqref{eq:3d_dynamics}, and~\eqref{eq:3d_constraints}), and uses additional arguments (e.g., cycle argument) to find the structure of optimal solutions. Thus, the conclusions in~\cite{sussmann1995} regarding maximals of minimum time problem will be leveraged throughout the following discussions to constuct $E^*(t_f)$. However, we note that the other arguments that are not based on Eqs.~\eqref{eq:thm1-1} and~\eqref{eq:thm1-2}, such as cycle argument, cannot be used for our reachability analysis.
%
%
\subsection{Reachability of 3D curves with prescribed curvature bound}
In this subsection, we present the reachability analysis of 3D curves with a prescribed curvature bound. The analysis is organized into two parts: first, we provide an implicit description of $\mathcal{G}_3(t_f)$ by proving that every boundary point can be reached by curves belonging to the classes CCC, CSC, their subsegments, or H. Second, we present the complete description of $\pi_3\left(\mathcal{G}_3(t_f)\right)$.
\subsubsection{Implicit description of the reachability set}
\begin{thm} \label{thm:6D}
	Any boundary point of $\mathcal{G}_3(t_f)$ can be reached by the curves of length $t_f$ satisfying the followings: 
	\begin{enumerate}
		\item it belongs to one of the following classes: CCC, CSC, their subsegments, or H, 
		\item if it is of class CCC, it is a planar curve, 
		\item if it is of class C, it is of length $\leq 4\pi$, 
		\item if it is of classes CS, SC, or CC, the C components have lengths $\leq 2\pi$, 
		\item if it is of classes CSC or CCC, the C components have lengths $< 2\pi$.
	\end{enumerate}
\end{thm}
\begin{pf*}{Proof.}
Proof of this theorem is largely an application of Theorem~\ref{thm:1} to the system Eq.~\eqref{eq:3d_dynamics}. We note that Eq.~\eqref{eq:3d_dynamics} satisfies Assumption~\ref{assum:01} because it is a smooth control affine system with compact and convex restraint set, $\mathbb{B}^3$. (Recall that the maps $\boldsymbol{w} \mapsto \boldsymbol{y} \times \boldsymbol{w}$ correspond to skew-symmetric linear transformations of $\mathbb{R}^3$.)

As introduced previously,~\cite{sussmann1995} provides a detailed description of the maximals of the minimum time problems described by Eqs.~\eqref{eq:3d_cost},~\eqref{eq:3d_dynamics}, and~\eqref{eq:3d_constraints}. In this regard, some portions of the proof of this theorem inherently have considerable overlap with the proof of Theorem~1 in~\cite{sussmann1995}. Therefore, during the proof of this theorem and the subsequent ones, the notations we define will follow the same ones used in the proof of Theorem~1 in~\cite{sussmann1995}. From the proof of Theorem~1 in~\cite{sussmann1995}, the following modifications will be made: (i) we extend the analysis to consider $\textbf{P}_{\boldsymbol{x}_f, \phi}$ problems with $\phi <\> 0$ as well, (ii) cycle argument should not be used, and (iii) we can exclude some points in $E^*(t_f)$ such that there exists an alternative trajectory reaching the same endpoint through a controller that is not maximal. The first and second items are related to construction of $E^*(t_f)$, and the last is to screen some maximal responses that their endpoints lie interior to $\mathcal{G}_3(t_f)$. The proof is organized as follows. First, in the `Basic setup' part, we will introduce the basic notations and summarize the developments from~\cite{sussmann1995} that will be used in the subsequent part. Since the three modifications previously mentioned are irrelevant with this step, this part is essentially contained in~\cite{sussmann1995}. In the subsequent Part 2, we prove items (1) and (2) of this theorem with the three modifications. Proof of the items (3), (4), and (5) are presented in Part 3.

Part 1. Basic setup and the structure of `nice' trajectories

Let us begin by defining the Hamiltonian as follows: 
\begin{equation}
	\mathcal{H}(\boldsymbol{x}, \boldsymbol{y}, \boldsymbol{\lambda}, \boldsymbol{\mu}, \nu, \boldsymbol{w}) = \langle \boldsymbol{\lambda}, \boldsymbol{y} \rangle + \langle \boldsymbol{\mu}, \boldsymbol{y} \times \boldsymbol{w} \rangle + \nu.
\end{equation}
While the dynamics of the state variables $(\boldsymbol{x}, \boldsymbol{y})$ is governed by Eq.~\eqref{eq:3d_dynamics}, the evolution of the costate variables $(\boldsymbol{\lambda}, \boldsymbol{\mu}, \nu)$ are governed by the costate differential equations 
\begin{equation} \label{eq:costate}
	\dot{\boldsymbol{\lambda}} = \boldsymbol{0}, \enspace \dot{\boldsymbol{\mu}} = -\boldsymbol{\lambda} - \boldsymbol{w} \times \boldsymbol{\mu}, \enspace \dot{\nu} = 0. 
\end{equation}
In terms of the notations used in the previous sections, $(\boldsymbol{\lambda}, \boldsymbol{\mu})$ corresponds to the adjoint variable $-\boldsymbol{p}$ and $\nu$ corresponds to $-p_0 \phi$. From here on, we will denote $\boldsymbol{\Xi} = (\boldsymbol{x}, \boldsymbol{y}, \boldsymbol{\lambda}, \boldsymbol{\mu}, \nu, \boldsymbol{w})$ satisfying Eqs.~\eqref{eq:3d_dynamics} and~\eqref{eq:costate} by \textit{lifted controlled trajectories}, and \textit{lifted controlled arcs} if it is defined on a compact interval. A lifted controlled trajectory is said to be \textit{$\mathcal{H}$-minimizing} if $\mathcal{H}\left( \boldsymbol{\Xi}(t) \right) = \min \{ \mathcal{H}(\boldsymbol{x}, \boldsymbol{y}, \boldsymbol{\lambda}, \boldsymbol{\mu}, \nu, \boldsymbol{v}) : \boldsymbol{v} \in \mathbb{B}^3 \}$ for a.e. $t$. 

Then a trajectory $(\boldsymbol{x}(t), \boldsymbol{y}(t))$ resulting from a control $\boldsymbol{w}(t)$ defined on a compact interval $[0, t_f]$ is said to satisfy the PMP conditions if there exists a lifted controlled arc $\boldsymbol{\Xi} = (\boldsymbol{x}, \boldsymbol{y}, \boldsymbol{\lambda}, \boldsymbol{\mu}, \nu, \boldsymbol{w})$ such that: 
\begin{enumerate}
	\item satisfies costate differential equation Eq.~\eqref{eq:costate},
	\item $\boldsymbol{\Xi}$ is $\mathcal{H}$-minimizing, 
	\item $\mathcal{H}\left( \boldsymbol{\Xi} \right) = 0$ for a.e. $t$, 
	\item $(\boldsymbol{\lambda}, \boldsymbol{\mu}) \neq \boldsymbol{0}$ as a covector on $M$ for all $t$. 
\end{enumerate}
These are merely a restatement of the PMP conditions of $\textbf{P}_{\boldsymbol{x}_f, \phi}$ with arbitrary $\boldsymbol{x}_f$ and $\phi \neq 0$. Indeed, when studying the minimum time problems, or $\textbf{P}_{\boldsymbol{x}_f, \phi}$ problems with $\phi < 0$, in~\cite{sussmann1995}, $\nu \geq 0$ condition was present. This condition is omitted here because we consider $\phi > 0$ case as well.

The costate vector $\boldsymbol{\mu}$ is associated with the state variable $\boldsymbol{y} \in \mathbb{S}^2$, and the linear functional $T_{\boldsymbol{y}}\mathbb{S}^2 \mapsto \mathbb{R}$ defined as $\boldsymbol{v} \mapsto \langle \boldsymbol{\mu}, \boldsymbol{v} \rangle$ becomes zero if and only if $\boldsymbol{\mu} \times \boldsymbol{y} = \boldsymbol{0}$. (Recall that $T_{\boldsymbol{y}}\mathbb{S}^2$ is $\boldsymbol{y}^{\perp}$.) Hence, $(\boldsymbol{\lambda}, \boldsymbol{\mu})$ does not vanish as a covector on $M$ if and only if 
\begin{equation} \label{eq:nontriviality}
	\|\boldsymbol{\lambda}\| + \|\boldsymbol{\mu}(t) \times \boldsymbol{y}(t)\| > 0 \quad \text{for all } t
\end{equation}
Dependency of $\boldsymbol{\lambda}$ on $t$ is ignored since Eq.~\eqref{eq:costate} implies that $\boldsymbol{\lambda}$ is constant in $t$. In~\cite{sussmann1995}, this formulation of the nontriviality condition\textemdash{}rather than the ``naïve'' version that requires $\left( \boldsymbol{\lambda}, \boldsymbol{\mu} \right) \neq 0$\textemdash{}played a crucial role. This is because the maximum principle with the latter nontriviality condition provided no useful information, as every trajectory was maximal~(Remark 3.2,~\cite{sussmann1995}). Under the nontriviality condition in Eq.~\eqref{eq:nontriviality}, the remaining part of the proof of this theorem demonstrates how maximum principle provides nontrivial conditions to sort out the maximals.

By definition of cross product, $\langle \boldsymbol{\mu}, \boldsymbol{y} \times \boldsymbol{w} \rangle = \det\left( \boldsymbol{y}, \boldsymbol{w}, \boldsymbol{\mu} \right) = \det\left( \boldsymbol{\mu}, \boldsymbol{y}, \boldsymbol{w} \right) = \langle \boldsymbol{w}, \boldsymbol{\mu} \times \boldsymbol{y} \rangle$ and consequently, Hamiltonian can be written as 
\begin{equation}
\mathcal{H}(\boldsymbol{\Xi}) = \langle \boldsymbol{\lambda}, \boldsymbol{y} \rangle + \langle \boldsymbol{w}, \boldsymbol{\mu} \times \boldsymbol{y} \rangle + \nu
\end{equation}
It is directforward that if $\boldsymbol{y} \times \boldsymbol{\mu} \neq \boldsymbol{0}$, $\mathcal{H}$ is minimized by $\boldsymbol{w} = -\frac{\boldsymbol{\mu} \times \boldsymbol{y}}{\|\boldsymbol{\mu} \times \boldsymbol{y}\|}$. If $\boldsymbol{y} \times \boldsymbol{\mu} = \boldsymbol{0}$, then the entire $\mathbb{B}^3$ becomes the set of minimizers. 
By defining the followings, 
\begin{equation}
	\varphi = \langle \boldsymbol{\lambda}, \boldsymbol{y} \rangle, \quad \boldsymbol{W} = \boldsymbol{y} \times \boldsymbol{\mu}
\end{equation} 
the minimization condition for the both cases can be expressed as 
\begin{equation}
	\|\boldsymbol{W}\| \boldsymbol{w} = \boldsymbol{W}
\end{equation}
Then the nontriviality condition and zero Hamiltonian condition can be equivalently written as 
\begin{equation} \label{eq:nontriviality_vanishing}
	\|\boldsymbol{\lambda}\| + \|\boldsymbol{W}\| > 0, \quad \|\boldsymbol{W}\| = \varphi + \nu
\end{equation}
In~\cite{sussmann1995}, a lifted controlled trajectory $\boldsymbol{\Xi}$ satisfying the PMP conditions is called \textit{nice} if $\boldsymbol{W}$ is everywhere nonvanishing. We will adopt the same terminology. A direct consequence of this definition is that nice trajectories have constant curvature of $1$ because $\boldsymbol{w} = -\frac{\boldsymbol{\mu} \times \boldsymbol{y}}{\|\boldsymbol{\mu} \times \boldsymbol{y}\|}$ implies $\|\dot{\boldsymbol{y}}\| = \|\boldsymbol{y} \times \boldsymbol{w}\| = 1$. 
Up to this point are the basic setup and statement of the conditions of PMP. 

After stating the conditions of PMP, proof of Theorem~1 in~\cite{sussmann1995} then first studied the structure of nice trajectories and subsequently analyzed how nice trajectories are connected up with trajectories such that $\boldsymbol{W} = 0$. The former was studied solely based on the four conditions of PMP we outlined and was irrelevant with the sign of $\nu$ and cycle argument. (Recall that $\nu \geq 0$ and cycle argument are not part of the PMP conditions here.) However, the latter was not. 

Our proof follows the similar step, but with the three aforementioned modifications. For the rest of part 1, we introduce the major conclusions\textemdash{}primarily on the structure of nice trajectories\textemdash{}made in~\cite{sussmann1995} that were solely grounded on the PMP conditions. As discussed in the previous subsection, these conclusions can be applied analogously to our proof. The summary is as follows: 
\begin{enumerate}[label=(\roman*)]
	\item if $\boldsymbol{\Xi}$ satisfies the PMP conditions, then $C \equiv \langle \boldsymbol{\lambda}, \boldsymbol{W} \rangle$ is conserved along $\boldsymbol{\Xi}$, 
	\item if $\boldsymbol{\Xi}$ is a nice lifted controlled trajectory, then 
	\begin{equation} \label{eq:torsion}
		\ddot{\boldsymbol{y}} + \boldsymbol{y} = \tau \boldsymbol{y} \times \dot{\boldsymbol{y}}
	\end{equation}
	where $\tau = -\frac{C}{(\varphi + \nu)^2}$, and $\tau$ can be calculated by solving the following ODE 
	\begin{equation} \label{eq:phiODE}
		\ddot{\varphi} + \varphi = \frac{C^2}{(\varphi + \nu)^3}
	\end{equation}
	with $\varphi > -\nu$ condition, 
	\item solution of Eq.~\eqref{eq:phiODE} with $\varphi > -\nu$ on a specified interval gives rise to a nice lifted controlled trajectory through Eq.~\eqref{eq:torsion}, 
	\item if $C \neq 0$, then solution of Eq.~\eqref{eq:phiODE} such that $\varphi > -\nu$ exists globally. 
\end{enumerate}
Item (i) above is about the entire lifted controlled trajectory $\boldsymbol{\Xi}$ while the others are about nice trajectories. These four items altogether enable the analysis of the structures of nice trajectories. In items (ii) and (iii), $\varphi > -\nu$ condition is essential because we are assuming a nice trajectory, so $\|\boldsymbol{W}\| = \varphi + \nu > 0$. From here on, unless with explicit mention, solution of Eq.~\eqref{eq:phiODE} will refer to those satisfying $\varphi > -\nu$ condition. In this context, a \textit{maximal solution} of Eq.~\eqref{eq:phiODE} refers to solutions that cannot be extended where $\varphi > -\nu$. Nice trajectories rising from a maximal solution of Eq.~\eqref{eq:phiODE} are referred to \textit{maximal nice trajectories}. In Eq.~\eqref{eq:torsion}, $\tau$ is the \textit{torsion} of the curve $t \mapsto \boldsymbol{x}(t)$ by definition. If a solution of Eq.~\eqref{eq:phiODE} exists on a specified interval, $\tau$ can be calculated through Eq.~\eqref{eq:phiODE} and the Fundamental Theorem of Curves implies that a smooth curve $t \mapsto \boldsymbol{x}(t)$ can be uniquely determined. (Recall that nice trajectories have constant curvature of $1$.) This gives rise to a nice trajectory. Thus, the subsequent analysis involves solutions of Eq.~\eqref{eq:phiODE} where $C \neq 0$ and $C = 0$ cases are handled separately. 

If $C \neq 0$, Eq.~\eqref{eq:phiODE} can be modified in terms of $\tau$, which becomes exactly Eq.~\eqref{eq:torsionODE} with $\zeta = \frac{2\nu}{\sqrt{|C|}}$. Since nice trajectories are globally defined in this case, there is no   need of further analysis how these connect up with trajectories such that $\boldsymbol{W} = 0$. Additional remark on $C \neq 0$ case is that the condition of $\nu \geq 0$ was present in~\cite{sussmann1995} because $\textbf{P}_{\boldsymbol{x}_f, \phi}$ problems with $\phi < 0$ were considered. Hence, the `H' curves were confined to the ones with $\zeta \geq 0$ but here we allow $\zeta < 0$ case as well. 

If $C = 0$, the general solution of Eq.~\eqref{eq:phiODE} is in the form of $\varphi(t) = A\cos(t - t_0)$, where $A$ and $t_0$ may vary over nice trajectories. If $|A| \geq |\nu|$, $\|\boldsymbol{W}\| = \varphi + \nu$ can indeed approach to $0$ and hence, nice trajectories can only exist on intervals of form $(t_0  - \delta, t_0  + \delta)$. If $|A| < \nu$, nice trajectory exists globally and does not exist if $|A| < -\nu$. Eq.~\eqref{eq:torsion} implies that this nice trajectory is a circular arc of curvature $1$. Additionally, $\{ \boldsymbol{y}, \boldsymbol{w}, \boldsymbol{y} \times \boldsymbol{w} \}$ forms an orthonormal basis of $\mathbb{R}^3$ on a nice trajectory because $\boldsymbol{y} \perp \boldsymbol{w}$ and $\|\boldsymbol{w}\| = 1$. But $C = \langle \boldsymbol{\lambda}, \boldsymbol{W} \rangle = \langle \boldsymbol{\lambda}, \|\boldsymbol{W}\| \boldsymbol{w} \rangle = 0$ and $\boldsymbol{W} \neq 0$ implies $\langle \boldsymbol{\lambda}, \boldsymbol{w} \rangle = 0$. Therefore, in terms of orthonormal basis, $\boldsymbol{\lambda} = \langle \boldsymbol{\lambda}, \boldsymbol{y} \rangle \boldsymbol{y} + \langle \boldsymbol{\lambda}, \boldsymbol{y} \times \boldsymbol{w} \rangle \boldsymbol{y} \times \boldsymbol{w} = \varphi \boldsymbol{y} + \dot{\varphi} \dot{\boldsymbol{y}}$ and $\|\boldsymbol{\lambda}\| = \varphi^2 + \dot{\varphi}^2 = A^2$. Thus, every nice trajectories of $C = 0$ rise from a solution $\varphi(t) = \|\boldsymbol{\lambda}\|\cos(t-t_0)$ of Eq.~\eqref{eq:phiODE}, where $t_0$ may vary over nice trajectories. 

Part 2. Proof of items (1) and (2)

In this part, we study how nice trajectories connect up with trajectories such that $\boldsymbol{W} = 0$ when it is not globally defined. This is essentially a process of identifying $S^*(t_f)$ and $E^*(t_f)$, but we alongside exclude some trajectories that have an alternative trajectory reaching the same endpoint through a controller that is not maximal.

Let us commence by considering a lifted controlled arc $\boldsymbol{\Xi}$ defined on an interval $I = [0, t_f]$ satisfying the PMP conditions. We define $L = \{ t \in I : \boldsymbol{W}(t) \neq 0 \}$ and its complement in $I$, $Q = I \setminus L = \{ t \in I : \boldsymbol{W}(t) = 0 \}$. Since $\|\boldsymbol{W}\| = \|\boldsymbol{y} \times \boldsymbol{\mu}\|$ is continuous, its inverse image of the open set $\mathbb{R} \setminus \{0\}$, $L$, is open in $I$.

Under these definitions, we study some characteristics of $\boldsymbol{W}$, beginning from its differentiability. 
Absolute continuity of $\boldsymbol{W} = \boldsymbol{y} \times \boldsymbol{\mu}$ follows from absolute continuity of $\boldsymbol{y}$ and $\boldsymbol{\mu}$ on compact interval $I$. Moreover, the following differentiation holds almost everywhere. 
\begin{equation} \label{eq:wdot}
	\begin{split}
		\dot{\boldsymbol{W}} & = \dot{\boldsymbol{y}} \times \boldsymbol{\mu} + \boldsymbol{y} \times \dot{\boldsymbol{\mu}} \\
		& = (\boldsymbol{y} \times \boldsymbol{w}) \times \boldsymbol{\mu} + \boldsymbol{y} \times (-\boldsymbol{\lambda} - \boldsymbol{w} \times \boldsymbol{\mu}) \\
		& = -\boldsymbol{y} \times \boldsymbol{\lambda} + (\boldsymbol{y} \times \boldsymbol{w}) \times \boldsymbol{\mu} - \boldsymbol{y} \times (\boldsymbol{w} \times \boldsymbol{\mu}) \\
		& = -\boldsymbol{y} \times \boldsymbol{\lambda}  + \boldsymbol{w} \times (\boldsymbol{\mu} \times \boldsymbol{y}) \\
		& = -\boldsymbol{y} \times \boldsymbol{\lambda} 
	\end{split}
\end{equation}
The last equality is because $\boldsymbol{w} = -\frac{\boldsymbol{\mu} \times \boldsymbol{y}}{\|\boldsymbol{\mu} \times \boldsymbol{y}\|}$ if $\boldsymbol{\mu} \times \boldsymbol{y} \neq 0$. 
Then $\boldsymbol{W}$ is an absolutely continuous function such that its derivative exists almost everywhere, which equals to a continuous function $-\boldsymbol{y} \times \boldsymbol{\lambda}$. This implies that $\boldsymbol{W}$ is everywhere differentiable and $\dot{\boldsymbol{W}} = -\boldsymbol{y} \times \boldsymbol{\lambda}$. 

For some $t_1 \in Q$, we have $\dot{\boldsymbol{W}}(t_1) = \boldsymbol{\lambda} \times \boldsymbol{y}(t_1)$. Consequently, 
\begin{equation} \label{eq:wnormdot}
\begin{split}
	\|\dot{\boldsymbol{W}}(t_1)\|^2 & = \|\boldsymbol{\lambda} \times \boldsymbol{y}(t_1)\|^2 \\
	& = \|\boldsymbol{\lambda}\|^2 \|\boldsymbol{y}(t_1)\|^2 - \langle \boldsymbol{\lambda}, \boldsymbol{y}(t_1) \rangle^2 \\
	& = \|\boldsymbol{\lambda}\|^2 - \varphi^2 \\
	& = \|\boldsymbol{\lambda}\|^2 - \nu^2 \\
\end{split}
\end{equation}
where Eq.~\eqref{eq:nontriviality_vanishing} and $\boldsymbol{W}(t_1) = 0$ are used in the last equality. Since $\boldsymbol{\lambda}$ and $\nu$ are constants, $\dot{\boldsymbol{W}}(t_1)$ are identical for all $t_1 \in Q$. 
This enables further classification of the zeros of $\boldsymbol{W}$. 

First, if $\|\boldsymbol{\lambda}\| < |\nu|$, then Eq.~\eqref{eq:wnormdot} is a contradiction. Thus, $\boldsymbol{W} \neq 0$ and $\boldsymbol{\Xi}$ is a single piece of nice arc if it exists. If $\|\boldsymbol{\lambda}\| > |\nu|$, then $\dot{\boldsymbol{W}}(t_1) \neq 0$ for all $t_1 \in Q$. So each $t_1$ is an isolated point of $Q$ and $Q$ is a discrete set. This implies that $\boldsymbol{\Xi}$ consists of concatenation of number of nice trajectories that are connected up with each other at points such that $\boldsymbol{W} = 0$. If $\|\boldsymbol{\lambda}\| = |\nu|$ then the zeros of $\boldsymbol{W}$ are not necessarily isolated and $\boldsymbol{\Xi}$ could consist of concatenations of nice trajectories with trajectories such that $\boldsymbol{W} \equiv 0$.

This classification motivates the subsequent analysis based on a case study which its structure is summarized in Fig.~\ref{fig01}. 
\begin{figure}[ht] 
\begin{center}
\resizebox{70mm}{!}{\includegraphics{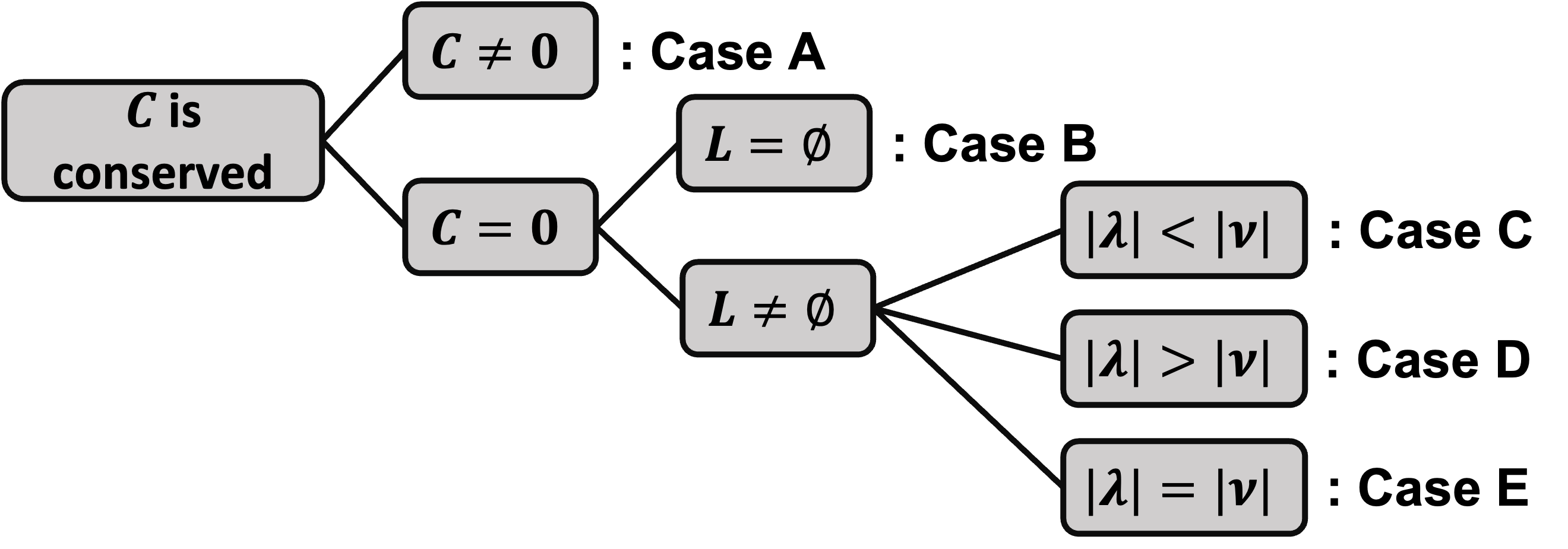}}
\caption{Structure of the case study analysis} \label{fig01}
\end{center}
\end{figure}
Case A, where $C \neq 0$, is already studied in part 1 concluding that such trajectories belong to class H. In this case, nice trajectories are globally defined and hence the entire $\boldsymbol{\Xi}$ is a nice arc. 
Case B corresponds to the case when $\boldsymbol{W} \equiv 0$. 
Case C is when nice trajectories can be globally defined, where the difference with case A is on $C = 0$ condition. 
Remaining cases D and E are relevant with how nice trajectories are connected up depending on $\boldsymbol{\lambda}$ and $\nu$, because they are not globally defined. 
Let us continue by addressing case B. 

Case B: $C = 0$ and $L = \emptyset$ 

This is when $\boldsymbol{W} \equiv 0$, which together with Eq.~\eqref{eq:nontriviality_vanishing} imply $\boldsymbol{\lambda} \neq 0$. Since $\boldsymbol{W} \equiv 0$, $\dot{\boldsymbol{W}} = -\boldsymbol{y} \times \boldsymbol{\lambda} = 0$ everywhere. But $\boldsymbol{y} \in \mathbb{S}^2$ and $\boldsymbol{\lambda}$ is nonzero, so $\boldsymbol{y} = \pm \frac{\boldsymbol{\lambda}}{\|\boldsymbol{\lambda}\|}$ and continuity of $\boldsymbol{y}$ implies that $\boldsymbol{y}$ is constant in time. Therefore, $\boldsymbol{\Xi}$ in case B corresponds to a straight segment, S. 

Case C: $C = 0$, $L \neq \emptyset$, and $\|\boldsymbol{\lambda}\| < |\nu|$ 

In this case, $\varphi(t) = \|\boldsymbol{\lambda}\|\cos(t-t_0)$ implies that a solution satisfying $\varphi > -\nu$ exists globally if $\|\boldsymbol{\lambda}\| < \nu$ and does not exist at all if $\|\boldsymbol{\lambda}\| < -\nu$. Hence, a nice trajectory globally exists if and only if $\|\boldsymbol{\lambda}\| < \nu$ and $\boldsymbol{\Xi}$ is a single piece of nice arc. It is already confirmed in part 1 that if $C = 0$, a single piece of nice trajectory corresponds to a circular arc of curvature $1$. It is noteworthy that the domain of $\boldsymbol{\Xi}$ need not be bounded because $\boldsymbol{W}$ never vanishes. 

Case D: $C = 0$, $L \neq \emptyset$, and $\|\boldsymbol{\lambda}\| > |\nu|$

Since $Q$ is a discrete set and each pieces of nice trajectories correspond to circular arcs of curvature 1, $\boldsymbol{\Xi}$ corresponds to a concatenation of multiple circular arcs. Suppose we denote $\boldsymbol{\Xi}_1 \boldsymbol{\Xi}_2 \dots \boldsymbol{\Xi}_{N}$ to represent such concatenations, where each $\boldsymbol{\Xi}_i$ are nice trajectories and $N$ is allowed to be infinite in principle. 

Since $\boldsymbol{W}$ is differentiable, Taylor's Theorem implies that for $t_1 \in Q$, $\boldsymbol{W}(t) = \dot{\boldsymbol{W}}(t_1)(t - t_1) + o(|t - t_1|)$ and $\boldsymbol{w}(t) = sgn(t - t_1)\frac{\dot{\boldsymbol{W}}(t_1)}{\|\dot{\boldsymbol{W}}(t_1)\|} + o(1)$ for $t \neq t_1$ near $t_1$. Therefore, the left and right side limits of $\boldsymbol{w}$ only change signs at $t_1$. This implies that the osculating plane of the two circular arcs being concatenated at $t_1$ coincides and their unit normal vectors are in the opposite directions. Hence, the state trajectory $\boldsymbol{x}(t)$ is a planar curve, concatenation of countable circular arcs. Moreover, $\boldsymbol{W} = 0$ at each point of concatenation, so $\boldsymbol{\Xi}_i$ for $2 \leq i \leq N-1$ correspond to maximal solutions of Eq.~\eqref{eq:phiODE}. Here, $\boldsymbol{\Xi}_1$ and $\boldsymbol{\Xi}_N$ are not necessarily maximal because $\boldsymbol{W}(0)$ and $\boldsymbol{W}(t_f)$ are not necessarily $0$.

From $\varphi(t) = \|\boldsymbol{\lambda}\|\cos(t-t_0)$, it is evident that maximal solutions are defined on nontrivial intervals of equal lengths $< 2\pi$. Consequently, $N$ has to be finite and the `middle' C curves corresponding to $2 \leq i \leq N-1$ have equal lengths $< 2\pi$, while the first and the last C's have lengths not greater than the middle C's. We will show after addressing the last case E that if $N > 3$, the endpoint $(\boldsymbol{x}(t_f), \boldsymbol{y}(t_f))$ lies interior to $\mathcal{G}_3(t_f)$. 

Case E: $C = 0$, $L \neq \emptyset$, and $\|\boldsymbol{\lambda}\| = |\nu|$ 

If $\nu = -\|\boldsymbol{\lambda}\|$, then it is easy to see that $\|\boldsymbol{W}\| = \|\boldsymbol{\lambda}\|\cos(t-t_0) + \nu$ can be nonnegative only at isolated points. Hence, no lifted controlled trajectory $\boldsymbol{\Xi}$ on a nontrivial interval can exist if $\nu = -\|\boldsymbol{\lambda}\|$. If $\nu = \|\boldsymbol{\lambda}\|$, then maximal solutions of Eq.~\eqref{eq:phiODE} are defined on open intervals of length $2\pi$. Since $L$ is open in $I$, it can be expressed as $L = \bigcup\limits_{J \in \mathcal{J}} J$ where $\mathcal{J}$ is a countable set of disjoint open subintervals of $I$. Then for any $J \in \mathcal{J}$ that does not contain $0$ and $t_f$, $\boldsymbol{W} = 0$ at its boundary points. Thus, such $J$ has a length of $2\pi$ and corresponds to the domain of a maximal nice trajectory. These maximal nice trajectories correspond to a single cycle. If the domain contains $0$ or $t_f$, then corresponding $\varphi$ does not have to be maximal, and naturally, such nice trajectory is a circular arc of length $\leq 2\pi$.

It is already shown in case B that the trajectories of $\boldsymbol{W} \equiv 0$ correspond to a straight line. Hence, the overall $\boldsymbol{\Xi}$ corresponds to a concatenation of multiple C and S components, where C components are full cycle (possibly multiple) if its domain does not contain $0$ or $t_f$. The cycles can be detached without perturbing the boundary conditions and then reattached tangentially to the initial point instead. Reattachment can be done in a way that the cycles and the first C component (if exists) altogether reduce to a single C component. This process reduces the overall trajectory to a curve of CSC or its subsegments, and consequently, points in $E^*(t_f)$ corresponding to case E can be reached by such curves. Unlike case D, the CSC curves need not be planar in general. 

Now, we will show that the endpoint $(\boldsymbol{x}(t_f), \boldsymbol{y}(t_f))$ lies interior to $\mathcal{G}_3(t_f)$ if $N > 3$ in case D. Suppose the corresponding $\boldsymbol{x}(t)$ of $\boldsymbol{\Xi}$ in case D is planar, consists of multiple C's with more than three components, denoted by $C_1 C_2 C_3 \dots C_{N-1} C_{N}$. According to the conclusion in case D, $C_2$ and $C_3$ have same lengths $< 2\pi$, and $C_1$ has a length not greater than the two. Thus, the sum of the lengths of $C_1$ and $C_3$ is larger than $C_2$ and none of them contains a cycle. Consequently, Lemma 2 in~\cite{patsko2003} asserts that there exists an alternative CCCCC curve, planar, such that its length, initial and terminal locations and directions coincide with $C_1 C_2 C_3$. Furthermore, it was proved that the lengths of the third and fourth C's are different in such CCCCC curves. By replacing the $C_1 C_2 C_3$ component by such curve, the resulting $CCCCC C_4 \dots C_{N-1} C_N$ curve has identical initial and terminal locations and directions with $C_1 C_2 C_3 \dots C_{N-1} C_{N}$. 
If this alternative trajectory is another outcome of case D, then the middle C's must have same lengths. However, this is not true. But along the conclusions of every other cases among A to E, we find that such trajectory is not a possible outcome of the PMP conditions. Therefore, the endpoints of the trajectories in case D with $N > 3$ have alternative trajectories reaching the same endpoint by a controller that is not maximal. Consequently, they cannot lie on the boundary of $\mathcal{G}_3(t_f)$. This completes the proof of items (1) and (2). Next, we prove the remaining items by identifying certain trajectories that posses an alternative trajectory reaching the same state endpoint $(\boldsymbol{x}(t_f), \boldsymbol{y}(t_f))$ by a controller that is not maximal.

Part 3. Proof of items (3), (4), and (5)

The stated classes in item (1) can be written out as H, S, C, CS, SC, CC, CSC, and CCC. Because some classes rise from multiple cases (e.g., CC trajectories can rise from cases D and E.), we will no longer follow the case study in Fig.~\ref{fig01} but follow the classification through different classes. 

First, consider a trajectory from class C with a length $> 4\pi$, where the corresponding $\boldsymbol{w}(t)$ takes a constant value $\boldsymbol{v}$ on $[0, t_f]$. This trajectory possesses at least two cycles. Then, one can find $t_1 < t_2 < t_3< t_4$ in $(0, t_f)$ such that $t_2 - t_1 = t_4 - t_3 = 2\pi$ and $t_3 - t_2 < 2\pi$. 
Subsequently, we can define an auxiliary control as follows. 
\begin{equation} \label{eq:aux_control}
	\tilde{\boldsymbol{w}}(t) = \begin{cases} 
		\boldsymbol{v}, \quad t \in [0, t_1) \\
		-\boldsymbol{v}, \quad t \in [t_1, t_2) \\
		\boldsymbol{v}, \quad t \in [t_2, t_3) \\
		-\boldsymbol{v}, \quad t \in [t_3, t_4) \\
		\boldsymbol{v}, \quad t \in [t_4, t_f]
	\end{cases}
\end{equation}
The auxiliary trajectory corresponding to this control input is depicted in Fig.~\ref{fig02} (a). This CCCCC trajectory has identical initial and terminal locations and directions with the original trajectory. Since it is of concatenation of multiple C components, it can only rise from cases D and E. However, the length of the C component corresponding to interval $[t_2, t_3]$ is less than $2\pi$ while the ones corresponding to $[t_1, t_2]$ and $[t_3, t_4]$ are $2\pi$. This trajectory cannot be an outcome of case D since it contains C components of length $2\pi$ in the middle. Moreover, it contradicts the conclusion of case E that all C components in the middle are cycles because $t_3 - t_2 < 2\pi$. Consequently, this auxiliary trajectory cannot rise from a nontrivial solution of PMP. This proves the item (3) of this theorem. 
\begin{figure}[htbp]
	\centering
	\subfigure[Auxiliary trajectory resulting from control input in Eq.~\eqref{eq:aux_control}]{\centering \includegraphics[scale=0.22]{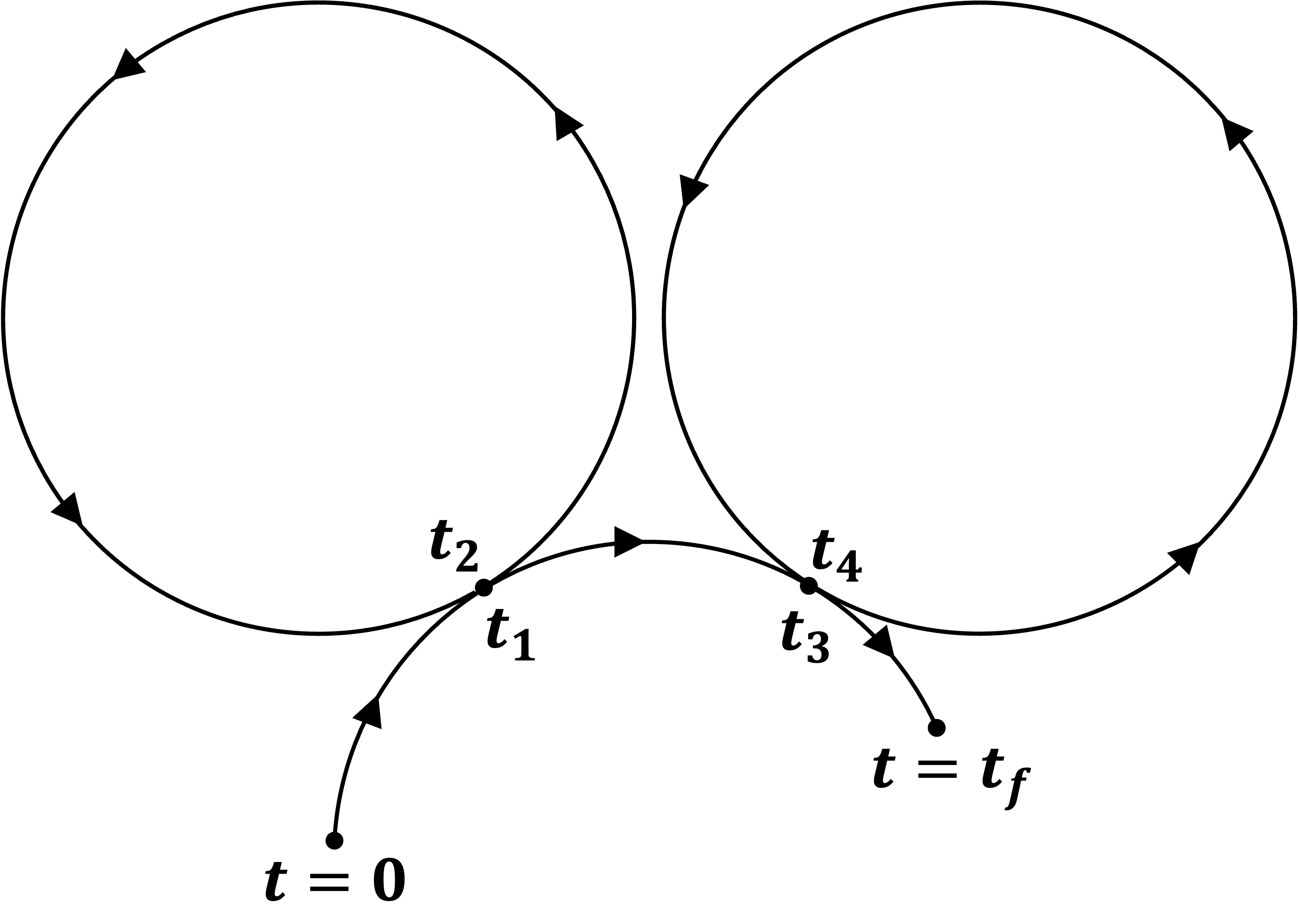}}
	\subfigure[Auxiliary trajectory of a CS curve]{\centering \includegraphics[scale=0.22]{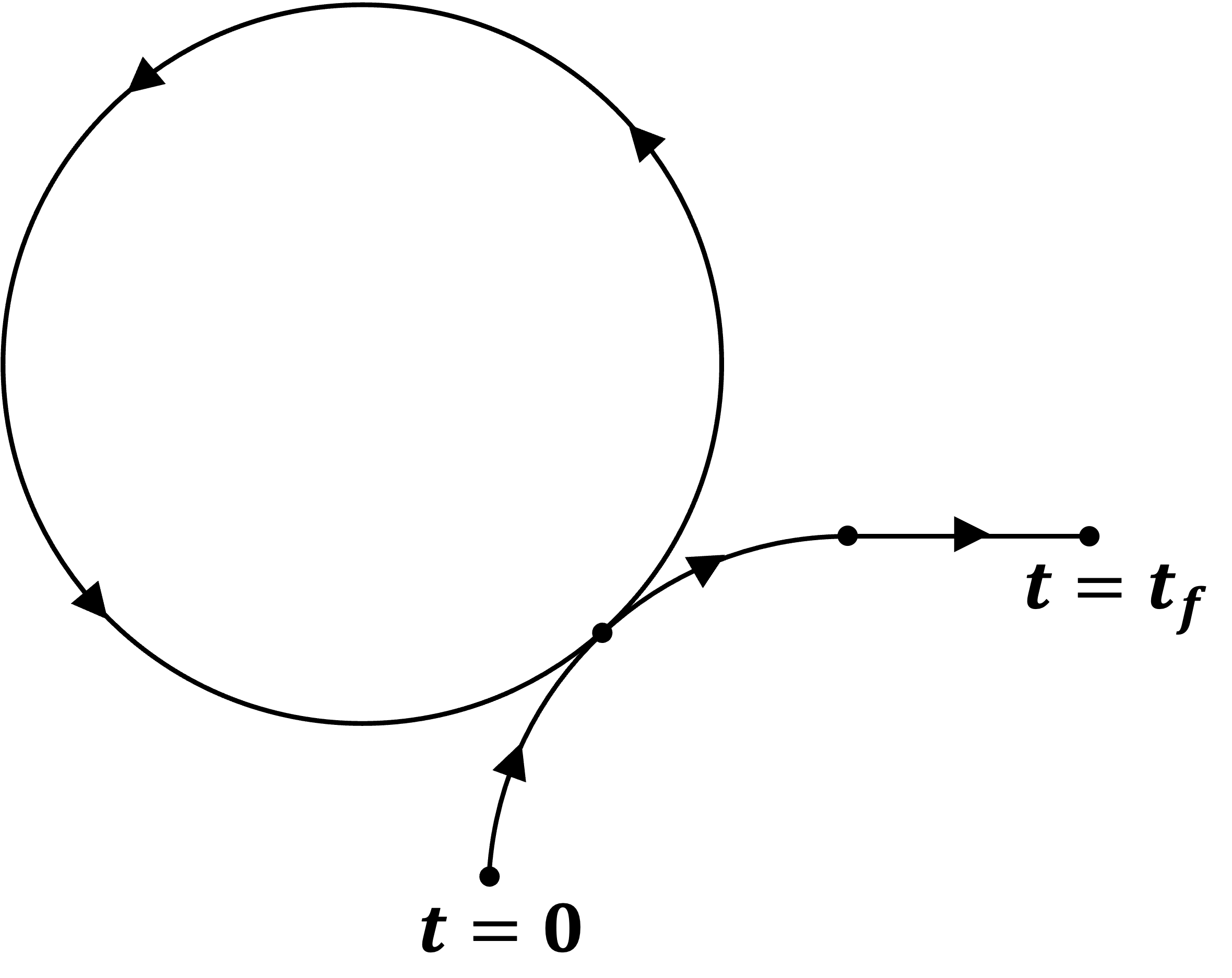}}
	\subfigure[Auxiliary trajectory of a CC curve]{\centering \includegraphics[scale=0.22]{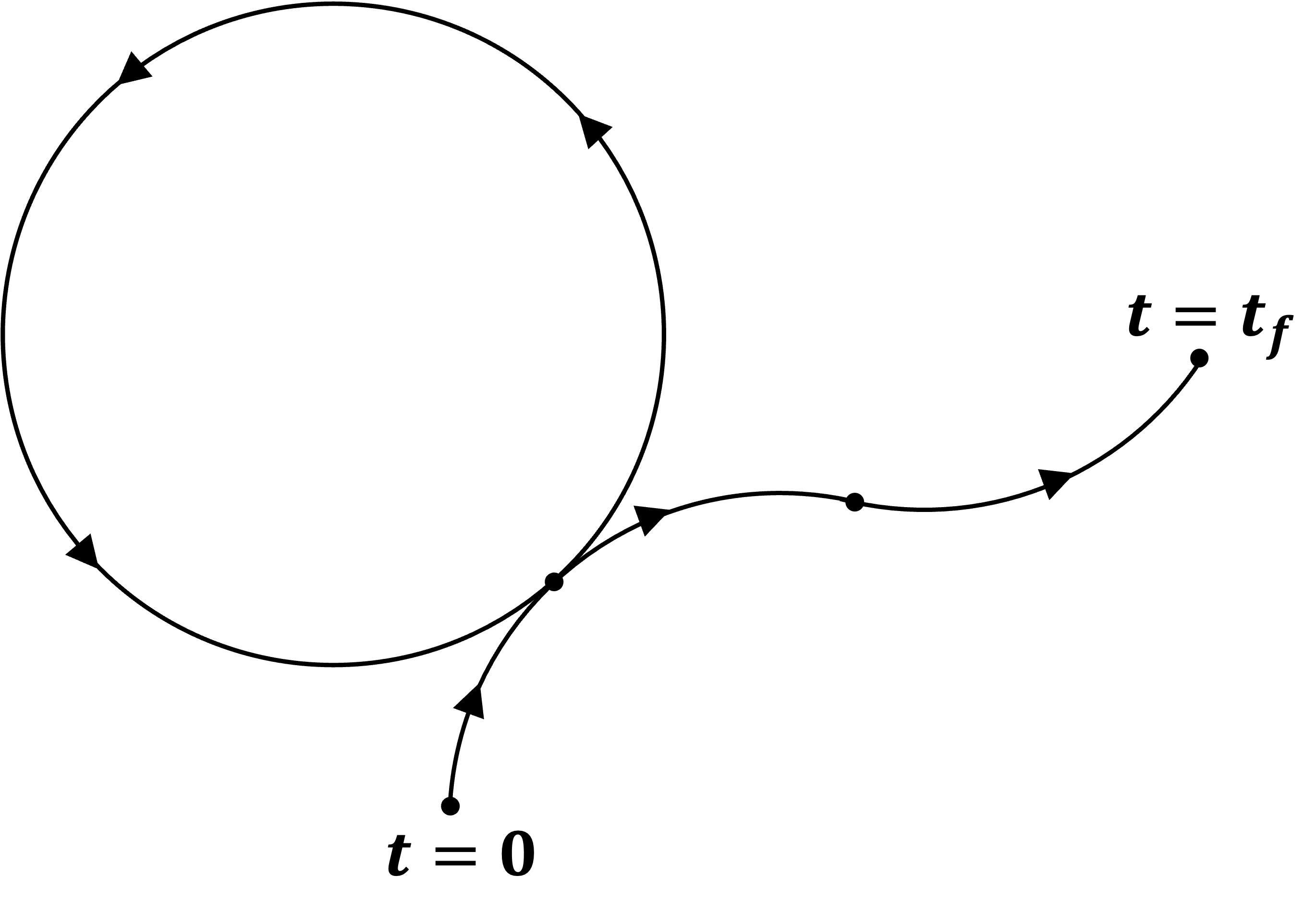}}
	\caption{Illustration of auxiliary trajectories that are generated from controllers that are not maximal. Among the state variables, only the $\boldsymbol{x}$ component is depicted.}
	\label{fig02}
\end{figure}

Next, suppose there is a trajectory belonging to class CS such that its C component has length $> 2\pi$. This implies that the C component possesses at least one cycle. Consequently, an auxiliary trajectory depicted in Fig.~\ref{fig02} (b) can be generated by similar construction as in Eq.~\eqref{eq:aux_control}, or by first detaching the cycle, and then reattaching it to the remaining C component. Since this CCCS trajectory contains an S segment, it can only rise from case E. But the third C component can be chosen to have length $< 2\pi$, which contradicts the conclusion of case E that all C components in the middle are cycles. Hence, such auxiliary trajectory cannot rise from a nontrivial solution of PMP. Analogous approach applies for SC trajectories as well. This proves that the endpoints of the trajectories belonging to the classes CS and SC lie interior to $\mathcal{G}_3(t_f)$ if the C component has length $> 2\pi$. 

Now, suppose there is a trajectory belonging to class CC which one of its C component has length $> 2\pi$. Similar to the previous cases, such C component possesses at least one cycle and an auxiliary trajectory belonging to class CCCC  depicted in Fig.~\ref{fig02} (c) can be generated by reattaching the cycle to either of the remaining C components. It should be noted that Fig.~\ref{fig02} (c) illustrates when the concatenation of C components is planar for elucidative purpose. However, it is not necessarily planar in general. (Recall that trajectories of class CC can also rise from case E.) Since this auxiliary trajectory contains a cycle, it can only rise from case E. But length of the third C component can be chosen to be less than $2\pi$, so that it cannot form a cycle. This contradicts the conclusion of case E. 

It remains to consider the CSC and CCC classes. If either of the C components of a CSC trajectory has length $\geq 2\pi$, it contains a cycle. Subsequently, an auxiliary trajectory can be generated by reattaching the cycle to one of the C components, as in the case of class CS in Fig.~\ref{fig02} (b). Then analogous statements prove that it cannot rise from a nontrivial solution of PMP. 

Lastly, we observe from part 2 that the CCC trajectories necessary for constructing the boundary of $\mathcal{G}_3(t_f)$ arise only from case D. This is because the CCC trajectories from case E can be reduced to CC and hence are redundant. Then, the conclusion of case D asserts that all C components have lengths $< 2\pi$. This completes the proof of the remaining items (3), (4), and (5). \qed
\end{pf*}
Items (1) and (2) of Theorem~\ref{thm:6D} collectively serve as the counterpart to Theorem~1 in~\cite{patsko2003} in 3D. Items (3), (4), and (5) of Theorem~\ref{thm:6D} are 3D counterparts of the developments in~\cite{buzikov2022} (Lemma 4 and Theorem 6). It should be noted that Theorem~1 in~\cite{patsko2003} and Theorem~\ref{thm:6D} in this paper, additional conditions to eliminate redundant trajectories with endpoints in the interior of reachability sets are required for complete determination of the boundary. For 2D curves,~\cite{patsko2022} provided such additional conditions and enabled complete determination of $\mathcal{G}_2(t_f)$. We leave the development of the counterpart of~\cite{patsko2022} for 3D curves as future work.

An additional remark is that the congruous relationship previously noted between~\cite{dubins1957} and~\cite{patsko2003}\textemdash{}both identifying classes CSC and CCC under 2D setting\textemdash{}is also observed in 3D settings. The structure of minimizers identified in \cite{sussmann1995}, which include CSC, CCC, their subsegments, and H, is also observed in Theorem~\ref{thm:6D}. However, the specific conditions under which these trajectories are optimal or lead to the boundary of the reachability set differs.

Acknowledging the subsequent studies (\cite{buzikov2022},~\cite{chen2023}, and~\cite{patsko2022}) that built upon the implicit description made in~\cite{patsko2003}, the 2D counterpart, Theorem~\ref{thm:6D} establishes a foundation for extending these studies into 3D settings. Moreover, this implicit description plays a crucial role in the subsequent discussion for complete determination of the position reachability set.
\subsubsection{Position reachability set}
\begin{figure}[htbp]
	\centering
	\subfigure{\centering \includegraphics[scale=0.18]{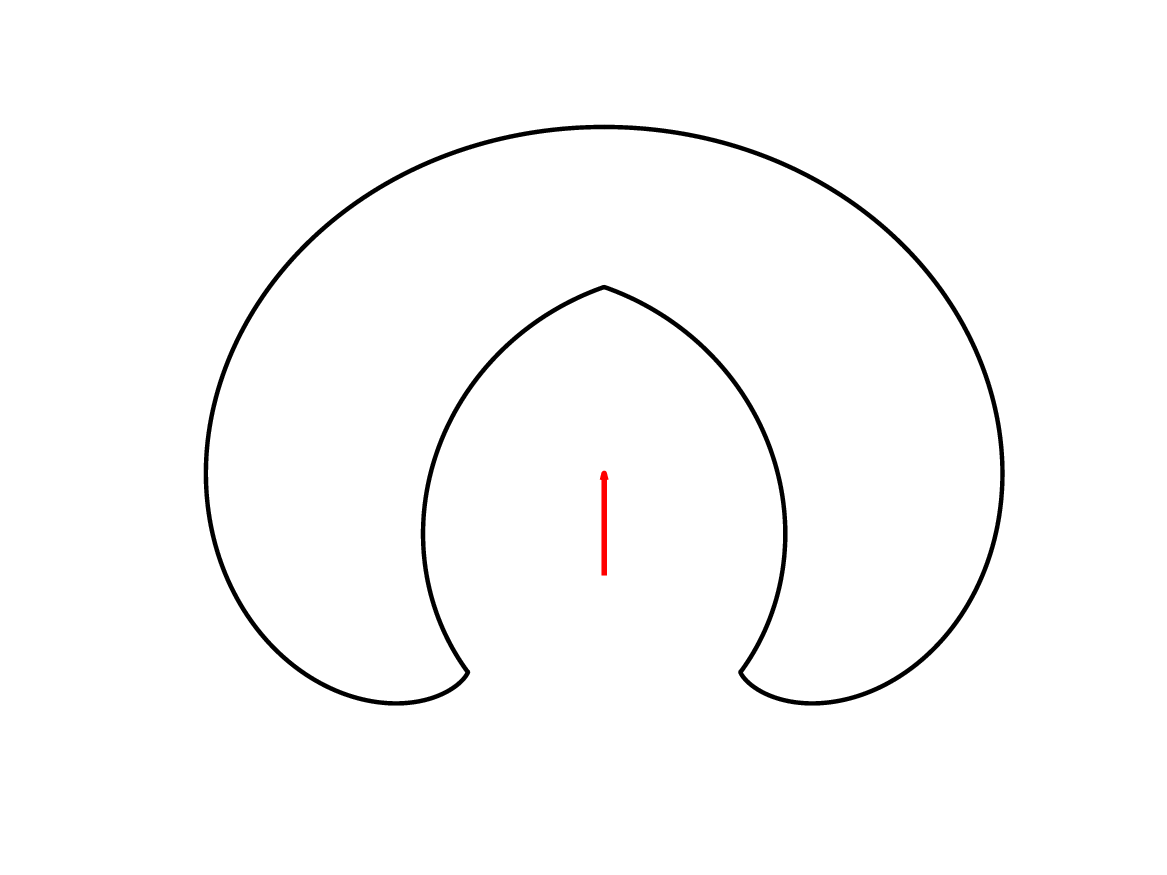}}
	\subfigure{\centering \includegraphics[scale=0.15]{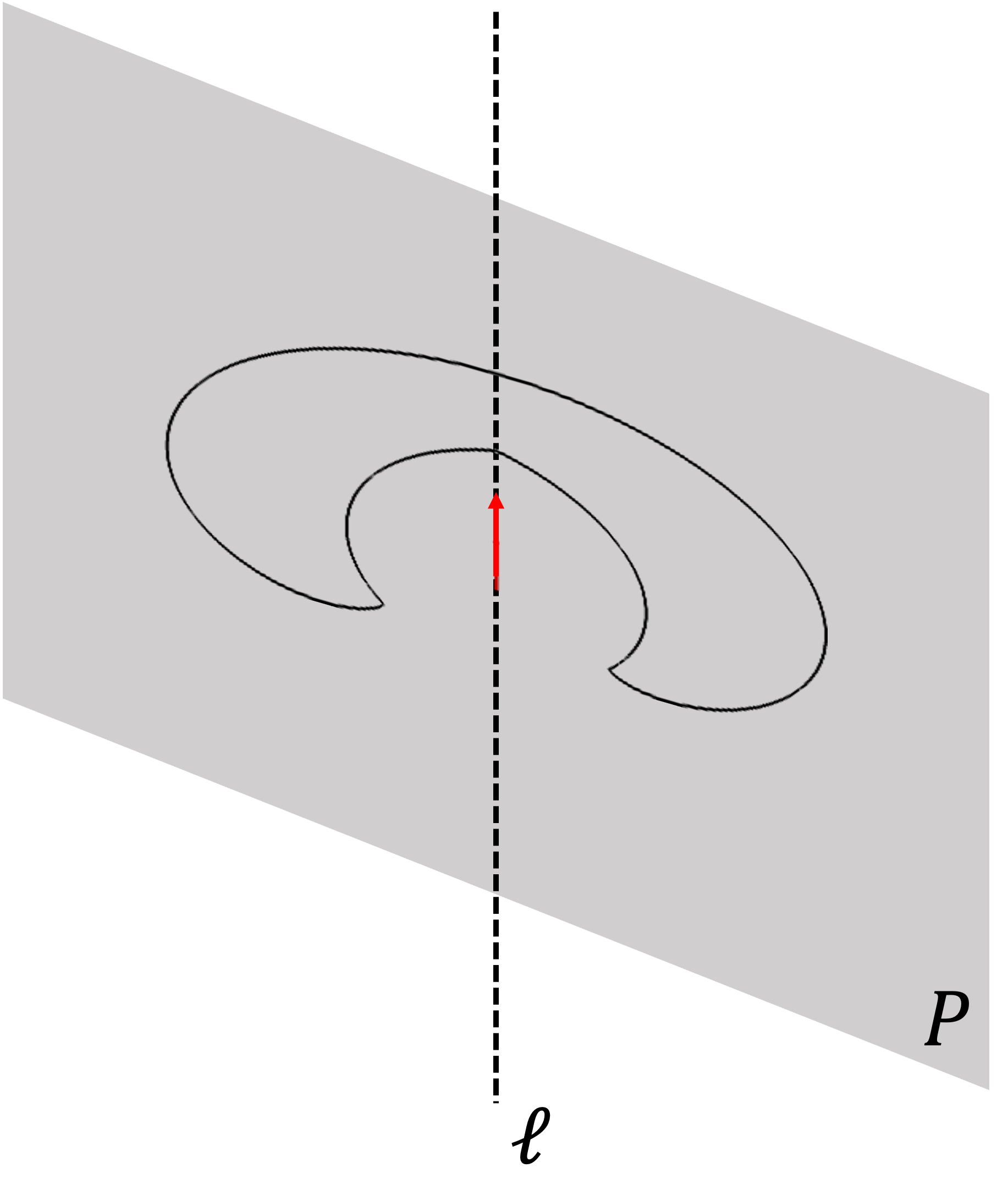}}
	\caption{Boundaries of $\pi_2\left(\mathcal{G}_2(t_f)\right)$ in $\mathbb{R}^2$ and in $P$. ($t_f = 1.4\pi$)}
	\label{fig2.5}
\end{figure}

Finally, we present the complete description of $\pi_3\left(\mathcal{G}_3(t_f)\right)$, or the position reachability set. The results expand the 2D case studied in~\cite{cockayne1975}. Unlike Theorem~\ref{thm:6D}, complete determination is done by additional geometric statements. We commence by stating the below theorem.
\begin{thm} \label{thm:3D}
Any boundary point of $\pi_3\left(\mathcal{G}_3(t_f)\right)$ can be reached by the curves of length $t_f$ belonging to the following classes: CS, CC, or their subsegments. Moreover, curves of CC that are necessary to reach the boundary points are limited to planar curves. 
\end{thm}
\begin{pf*}{Proof.}
The proof is based on application of Theorem~\ref{thm:2}. As stated below Theorem~\ref{thm:2}, the only additional condition compared to Theorem~\ref{thm:1} is that the adjoint response corresponding to the state variables being dropped by the projection vanishes as a covector at $t_f$. Under our setting, since $\pi_3: \mathbb{R}^3 \times \mathbb{S}^2 \mapsto \mathbb{R}^3$, we set $M_1 = \mathbb{R}^3$ and $M_2 = \mathbb{S}^2$ where $M_1$ and $M_2$ are as in Theorem~\ref{thm:2}. Then, we can continue upon the proof of Theorem~\ref{thm:6D} with an additional condition that $\boldsymbol{\mu}(t_f)$ is a zero covector in $T^*_{\boldsymbol{y}(t_f)}\mathbb{S}^2$. As confirmed along Eq.~\eqref{eq:nontriviality}, $\boldsymbol{\mu}(t_f)$ vanishes if and only if $\boldsymbol{W}(t_f) = 0$. Since the structure of nice trajectories are independent with this condition, part 1 in the proof of Theorem~\ref{thm:6D} holds analogously. Thus, it suffices to redo the proof of part 2 with an additional condition of $\boldsymbol{W}(t_f) = 0$. 
	
Along the case study outlined in Fig.~\ref{fig01}, application of $\boldsymbol{W}(t_f) = 0$ condition to cases A, B, and C follows trivially. Case A, which corresponds to the trajectories of class H, can be disregarded because $\boldsymbol{W}(t_f) = 0$ directly implies $C = \langle \boldsymbol{\lambda}, \boldsymbol{W} \rangle = 0$. Case B, which is when $\boldsymbol{W} \equiv 0$, essentially implies $\boldsymbol{W}(t_f) = 0$ and there is nothing to do more. Case C is disregarded because $\boldsymbol{W}$ cannot vanish in this case. Next, we need to address cases D and E. 

In case D, a lifted controlled arc $\boldsymbol{\Xi}$ satisfying the PMP conditions consists of concatenation of finite number of nice trajectories, denoted as $\boldsymbol{\Xi}_1 \boldsymbol{\Xi}_2 \dots \boldsymbol{\Xi}_{N}$. Moreover, $\boldsymbol{\Xi}_i$ for $2 \leq i \leq N-1$ correspond to a maximal solution of Eq.~\eqref{eq:phiODE} and therefore have the same lengths. But $\boldsymbol{W}(t_f) = 0$ further indicates that $\boldsymbol{\Xi}_N$ also corresponds to a maximal solution. Therefore, all $\boldsymbol{\Xi}_i$'s, except possibly $\boldsymbol{\Xi}_1$, must have the same lengths. Additionally, since the endpoints $(\boldsymbol{x}(t_f), \boldsymbol{y}(t_f))$ of trajectories in $N > 3$ cases lie interior to $\mathcal{G}_3(t_f)$, their projections by $\pi_3$ also lie interior to $\pi_3\left( \mathcal{G}_3(t_f) \right)$, given that projection is an open map. In other words, $N$ is at most 3 for $\boldsymbol{x}(t_f)$ to lie on the boundary of $\pi_3\left( \mathcal{G}_3(t_f) \right)$. If $N = 3$ and if we write $C_1 C_2 C_3$, the lengths of $C_2$ and $C_3$ are the same. Consequently, Lemma~2 in~\cite{patsko2003} implies that there exists an alternative CCCCC trajectory. During the proof of Theorem~\ref{thm:6D}, we confirmed that this cannot be an outcome of the PMP conditions. Hence, trajectories rising from case D required for construction of $\pi_3\left(\mathcal{G}_3(t_f)\right)$ are confined to C or planar CC trajectories such that the length of the second C is not less than the first.
	
In Case E of the proof of Theorem~\ref{thm:6D}, it is proved that each of the nice trajectories, corresponding to the C components, is of length $2\pi$ and rises from a maximal solution of Eq.~\eqref{eq:phiODE}, provided that its domain does not include the initial time $0$ or the final time $t_f$. However, if $\boldsymbol{W}(t_f) = 0$, the last C component corresponds to a maximal nice trajectory as well, and its length is $2\pi$. Consequently, the same process as in the proof of case E in Theorem~\ref{thm:6D} can be applied: such cycle can be reattached tangentially to the initial point. This reduces the overall curve into curves of CS or its subsegments. \qed
\end{pf*}
From here on, we denote by $\ell$ the line in $\mathbb{R}^3$ passing through the point $\boldsymbol{x}(0)$ with direction $\boldsymbol{y}(0)$. An arbitrary plane in $\mathbb{R}^3$ containing $\ell$ is indicated as $P$. Recall that $\pi_2\left( \mathcal{G}_2(t_f) \right)$ represents the set of all endpoints of curvature bounded paths in $\mathbb{R}^2$. Henceforth, under the same curvature bound setting, we abuse the notation and indicate by $\pi_2\left( \mathcal{G}_2(t_f) \right)$ the set of all endpoints of planar curves that lie entirely in $P$. In other words, the ambient space is considered as $P$ instead of $\mathbb{R}^2$. Schematic illustration of boundaries of $\pi_2\left( \mathcal{G}_2(t_f) \right)$ in $\mathbb{R}^2$ and in $P$ are depicted in Fig.~\ref{fig2.5}.
Similarly, let us denote by $E_1^P(t_f)$ the set of endpoints of curves that lie entirely in $P$ and belong to class CS, CC, or their subsegments. Then Proposition~\ref{prop:1} first asserts that $E_1^P(t_f)$ contains the boundary of $\pi_2\left(\mathcal{G}_2(t_f)\right)$. Subsequently, by denoting as $E'^P(t_f)$ the excluded points through the five additional items in Proposition~\ref{prop:1}, the statement is strengthened as follows: \textit{$E_1^P(t_f) \setminus E'^P(t_f)$ contains the boundary of $\pi_2\left(\mathcal{G}_2(t_f)\right)$}. 

In~\cite{cockayne1975} where the ambient space was $\mathbb{R}^2$, it was confirmed that there is a finite number of regions enclosed by $E_1(t_f) \setminus E'(t_f)$, which there is a unique way to select the enclosed regions so that their union varies continuously in time. This argument facilitated complete determination of $\pi_2\left( \mathcal{G}_2(t_f) \right)$. Analogous argument on $E_1^P(t_f) \setminus E'^P(t_f)$ allows unique determination of $\pi_2\left(\mathcal{G}_2(t_f)\right)$ in $P$. 

Now, let us denote by $E_2^P(t_f)$ the set of endpoints in $P$ of all planar curves in the classes CS, CC, and their subsegments. The difference with $E_1^P(t_f)$ is that the curves need not be entirely contained in $P$ but are only required to be planar. Since $\pi_3\left( \mathcal{G}_3(t_f) \right)$ is invariant under rotation around the axis $\ell$, it can be generated by rotating $\mathcal{S}(t_f) \equiv \pi_3\left( \mathcal{G}_3(t_f) \right) \cap P$ around $\ell$. It is easy to see that $\pi_2\left( \mathcal{G}_2(t_f) \right) \subseteq \mathcal{S}(t_f)$ because the responses of Eq.~\eqref{eq:3d_dynamics} include the planar curves that lie in $P$. Moreover, $\mathcal{S}(t_f)$ varies continuously in time if and only if $\pi_3\left( \mathcal{G}_3(t_f) \right)$ does.

\begin{figure}[htbp]
	\centering
	\subfigure[$t_f = 0.5\pi$]{\centering \includegraphics[scale=0.13]{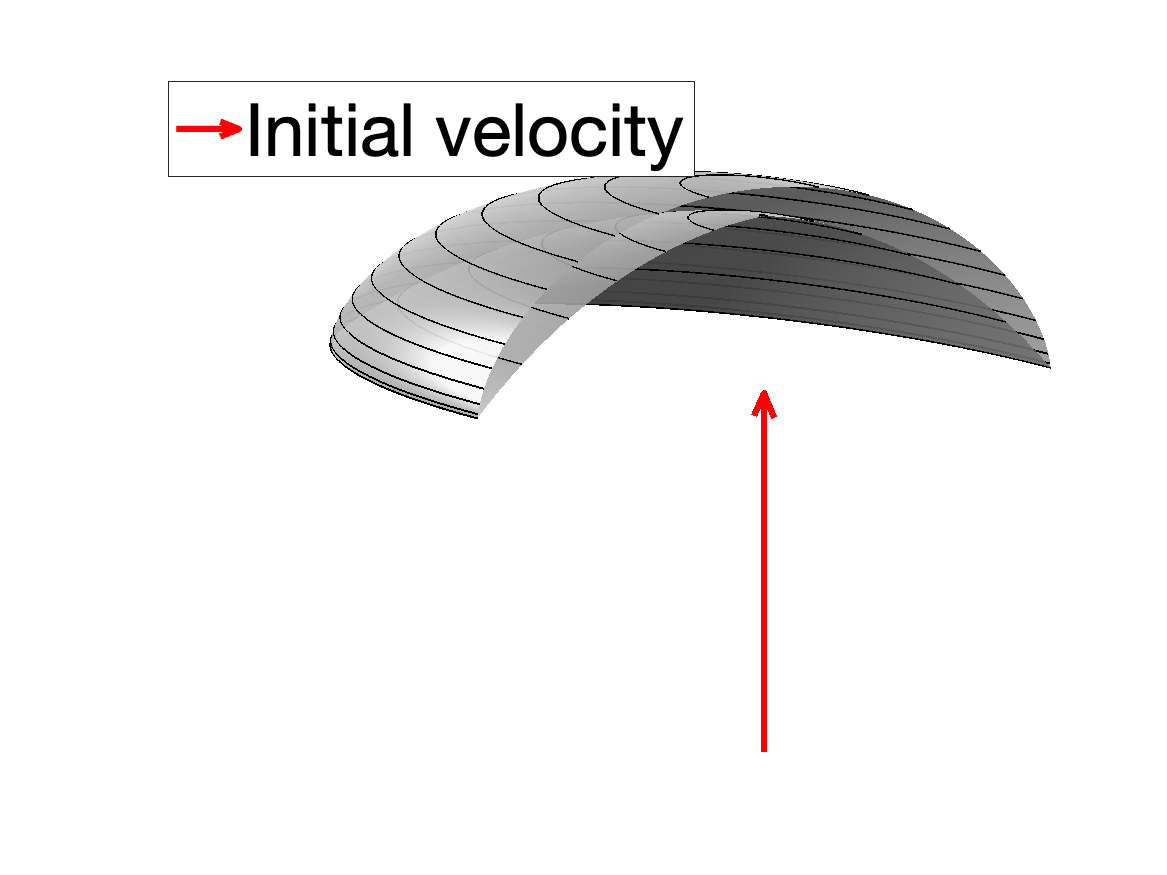}}
	\subfigure[$t_f = 0.8\pi$]{\centering \includegraphics[scale=0.13]{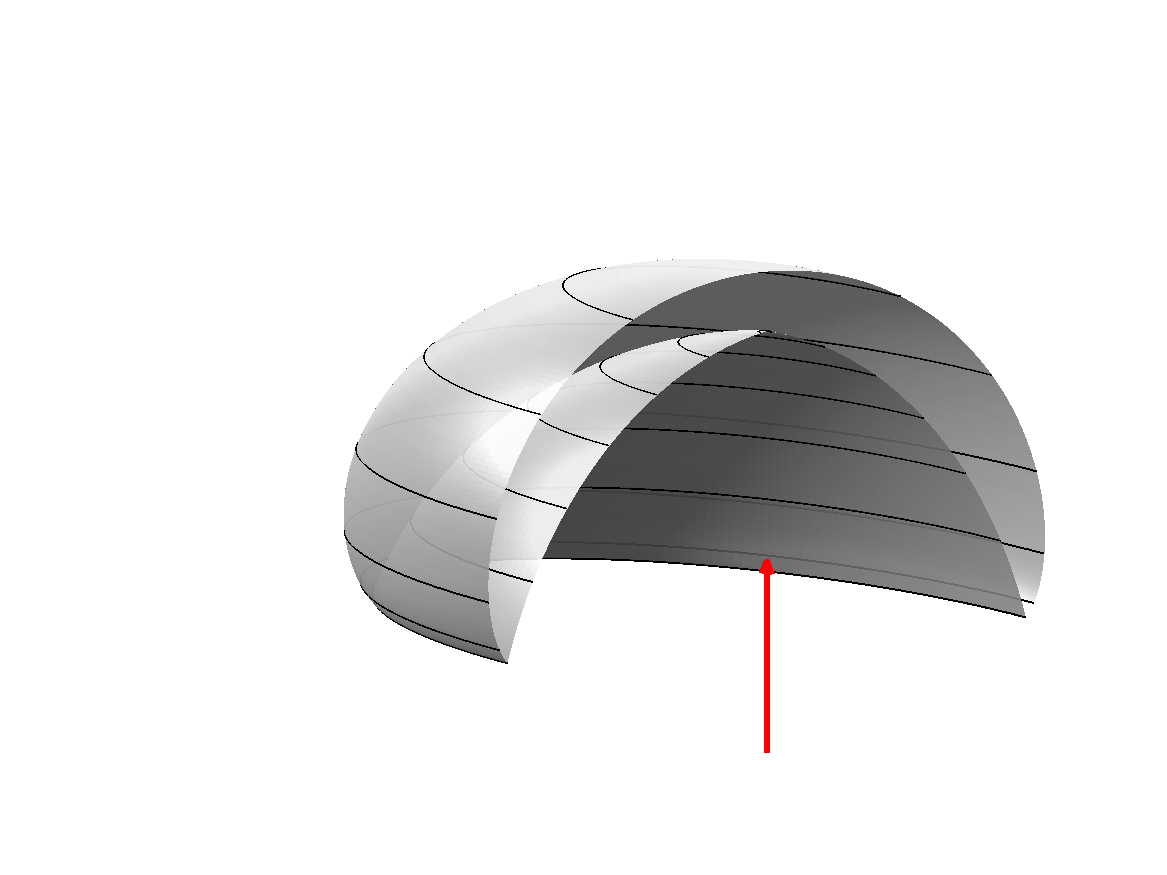}}
	\subfigure[$t_f = 1.5\pi$]{\centering \includegraphics[scale=0.13]{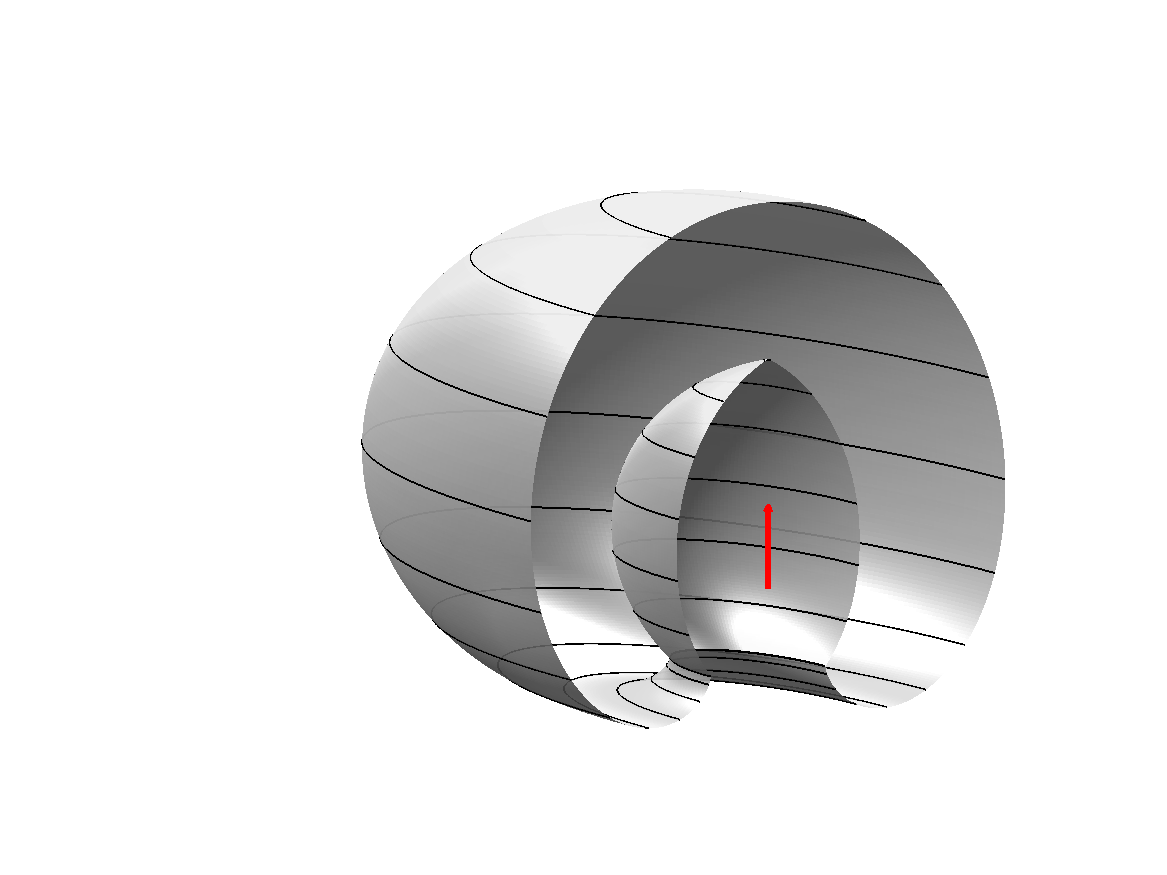}}
	\subfigure[$t_f = 1.8\pi$]{\centering \includegraphics[scale=0.13]{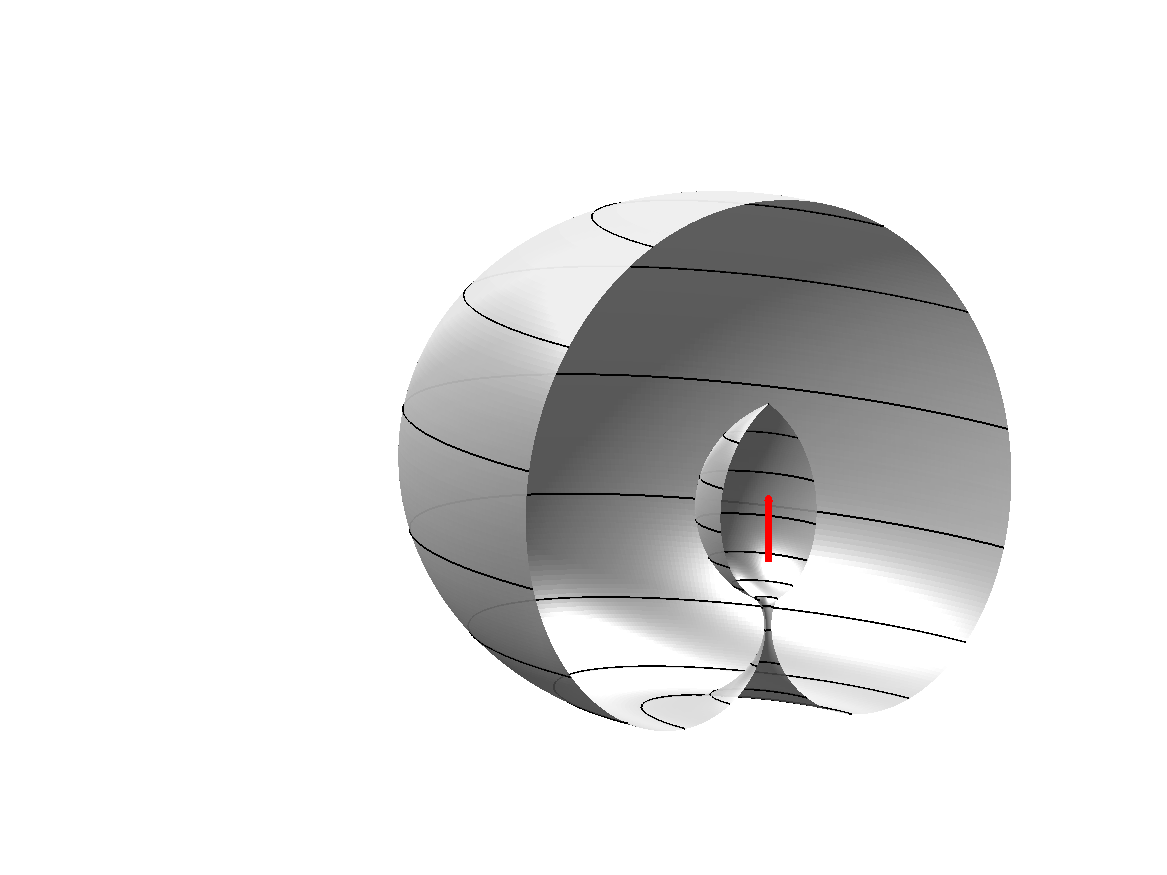}}
	\subfigure[$t_f = 1.9\pi$]{\centering \includegraphics[scale=0.13]{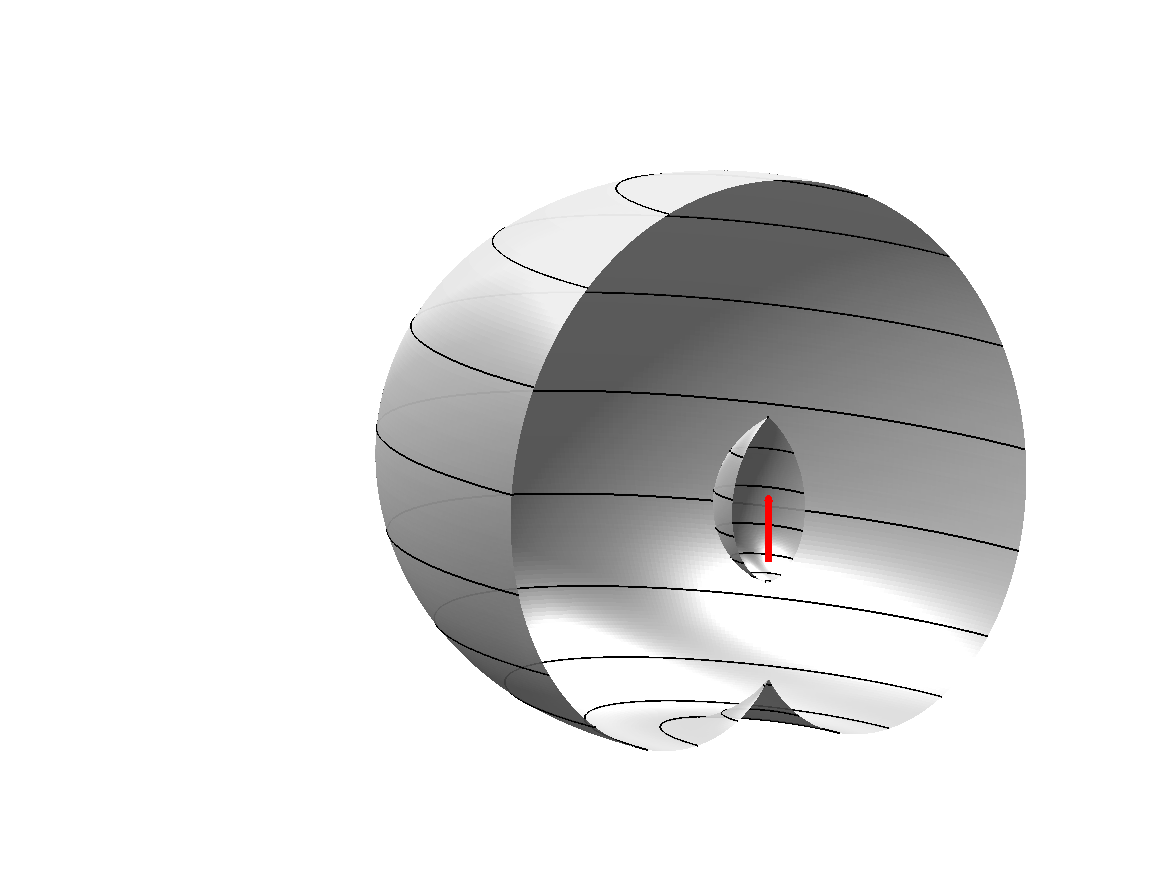}}
	\subfigure[$t_f = 3\pi$]{\centering \includegraphics[scale=0.13]{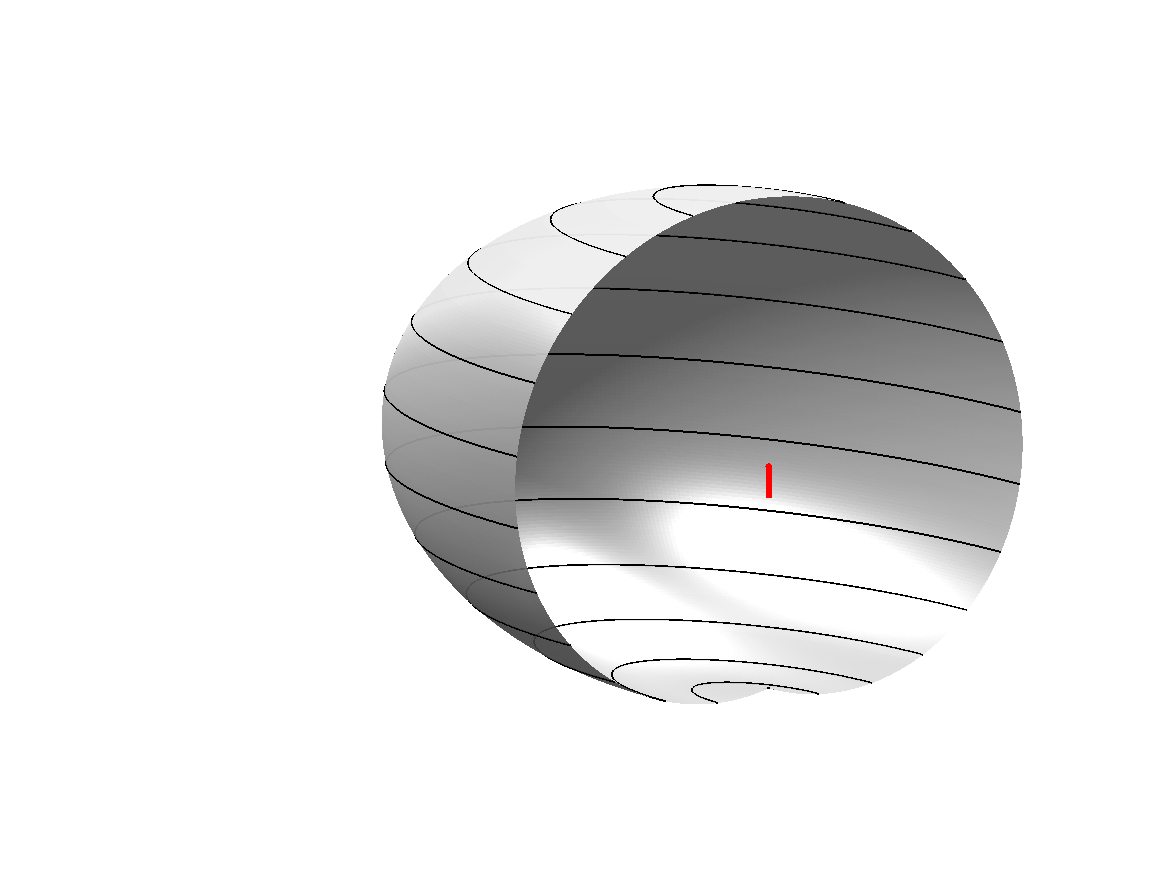}}
	\caption{Cross-sectional view of $\pi_3\left(\mathcal{G}_3(t_f)\right)$ over multiple $t_f$.}
	\label{fig03}
\end{figure}

Subsequently, we show that $\pi_2\left( \mathcal{G}_2(t_f) \right) = \mathcal{S}(t_f)$. It follows from Theorem~\ref{thm:3D} that $E_2^P(t_f)$ contains the boundary of $\mathcal{S}(t_f)$, and the union of $E_2^P(t_f)$ over all $P$ contains the boundary of $\pi_3\left( \mathcal{G}_3(t_f) \right)$. We introduce the following lemma.
\begin{lem}
	$E_1^P(t_f) = E_2^P(t_f)$. 
\end{lem}
\begin{pf*}{Proof.}
	From definition, it is evident that $E_1^P(t_f) \subseteq E_2^P(t_f)$. We will show the inclusion relationship in the converse direction. Suppose a point $\boldsymbol{z}$ in $E_2^P(t_f)$ is reached by a planar curve in classes CS, CC, or their subsegments that lies entirely in a plane $P'$. It follows from definition that $\boldsymbol{z}$ is in both $P$ and $P'$. Since $\ell$ is the tangent line of the curve at $t = 0$, it follows $\ell \subseteq P'$. The claim is trivial if $P = P'$, so we assume $P \neq P'$. Then $\ell$ is a line entirely contained in two distinct planes $P$ and $P'$. Thus, $\ell = P \cap P'$ and $\boldsymbol{z} \in \ell$. Then one can rotate the entire curve\textemdash{}which originally lies in $P'$\textemdash{}around $\ell$ until the curve lies in $P$. This maintains the endpoint on $\boldsymbol{z}$ because it is on the axis of rotation, $\ell$. The result ing curve is of length $t_f$, lies in $P$, belongs to class CS, CC, or their subsegments, and reaches $\boldsymbol{z}$. This proves that $\boldsymbol{z} \in E_1^P(t_f)$. \qed 
\end{pf*}
Now, the points in $E_1^P(t_f)$\textemdash{}or equivalently, $E_2^P(t_f)$\textemdash{}that violate any of the five items in Proposition~\ref{prop:1} lie interior to $\pi_2\left(\mathcal{G}_2(t_f)\right)$ and, consequently, interior to $\mathcal{S}(t_f)$ since $\pi_2\left( \mathcal{G}_2(t_f) \right) \subseteq \mathcal{S}(t_f)$. This strengthens Theorem~\ref{thm:3D} as follows: \textit{$E_2^P(t_f) \setminus E'^P(t_f)$ contains the boundary of $\mathcal{S}(t_f)$}. But given $E_1^P(t_f) \setminus E'^P(t_f)$, or equivalently $E_2^P(t_f) \setminus E'^P(t_f)$, there is a unique way to choose the enclosed regions in $P$ such that their union varies continuously in time. This union equals to $\pi_2\left( \mathcal{G}_2(t_f) \right)$ and $\mathcal{S}(t_f)$. This proves that $\mathcal{S}(t_f) = \pi_2\left( \mathcal{G}_2(t_f) \right)$ and, subsequently, the following theorem. This allows complete determination of $\pi_3\left(\mathcal{G}_3(t_f)\right)$ for all $t_f$ as in Fig.~\ref{fig03}.
\begin{thm} \label{thm:3D_complete}
	$\pi_3\left(\mathcal{G}_3(t_f)\right)$ is uniquely determined by rotating $\pi_2\left(\mathcal{G}_2(t_f)\right)$ around the initial velocity vector. 
\end{thm}
\begin{rem}
	(Remark on the general-dimensional case) The structure of maximal trajectories for the Markov-Dubins problem in general dimensions was studied in~\cite{Mittenhuber1998}, also based on PMP. The analysis extended beyond euclidean spaces to include noneuclidean spaces such as spheres, $\mathbb{S}^n$, and hyperbolic spaces, $\mathbb{H}^n$. It was proved that to solve the Markov-Dubins problem in $\mathbb{R}^n$ (resp. $\mathbb{S}^n$, $\mathbb{H}^n$), it suffices to solve it in $\mathbb{R}^3$ (resp. $\mathbb{S}^3$, $\mathbb{H}^3$). This suggests a promising avenue for applying the approach of this paper to derive conclusions in general dimensions. However, this requires to entirely reformulate the problem on the Lie group $SE(n)$ of Euclidean motions of $\mathbb{R}^n$. We leave this extension for future work.
\end{rem}
%
\section{Conclusion}\label{sec:04}
In this paper, we have presented the reachability analysis of three-dimensional curves with a prescribed curvature bound. Based on Pontryagin Maximum Principle and existing results on Markov-Dubins problem in dimension 3, we provided an implicit description of the boundary of the reachability set that any boundary point can be reached by curves belonging to the classes CSC, CCC, their subsegments, or H. We further proved that the position reachability set is a solid of revolution of its two-dimensional counterpart, the Dubins car, thereby enabling complete determination. These developments extend the existing literature on Dubins car into dimension 3.
%
\begin{ack}
We would like to extend our gratitude to Nearthlab for supports for this paper, and the anonymous reviewers for their insightful comments.
\end{ack}
\bibliographystyle{automatica}        
\bibliography{Ref_Papers.bib}
\end{document}